\documentclass{tran-l}

\usepackage{array,float}
\usepackage{amssymb}

\copyrightinfo{2006}{Benjamin Klopsch and Christopher Voll}

\newtheorem{thmx}{Theorem} 

\newtheorem{conjecture}[thmx]{Conjecture}

\newtheorem{theorem}{Theorem}[section] 

\newtheorem{proposition}[theorem]{Proposition} 
\newtheorem{lemma}[theorem]{Lemma}

\theoremstyle{definition}
\newtheorem{definition}[theorem]{Definition}
\newtheorem{example}[theorem]{Example}

\theoremstyle{remark}
\newtheorem{remark}[theorem]{Remark}

\numberwithin{equation}{section}

\newcommand{\F}{\ensuremath{\mathbb{F}}}
\newcommand{\N}{\ensuremath{\mathbb{N}}}

\newcommand{\bfb}{\ensuremath{\mathbf{b}}}
\newcommand{\bfl}{\ensuremath{\mathbf{l}}}

\newcommand{\bfB}{\ensuremath{\boldsymbol{B}}}
\newcommand{\bfF}{\ensuremath{\mathbf{F}}}
\newcommand{\bfR}{\ensuremath{\mathbf{R}}}
\newcommand{\bfU}{\ensuremath{\mathbf{U}}}
\newcommand{\bfX}{\ensuremath{\mathbf{X}}}

\newcommand{\calA}{\ensuremath{\mathcal{A}}}
\newcommand{\calH}{\ensuremath{\mathcal{H}}}
\newcommand{\calV}{\ensuremath{\mathcal{V}}}
\newcommand{\calW}{\ensuremath{\mathcal{W}}}

\newcommand{\upL}{\ensuremath{\textup{L}}}
\newcommand{\upR}{\ensuremath{\textup{R}}}

\newcommand{\disjunion}{\ensuremath{\; \dot{\cup} \;}}

\newcommand{\e}{\ensuremath{\gamma}}
\newcommand{\sym}{\ensuremath{\mathcal{S}}}

\DeclareMathOperator{\cha}{char}
\DeclareMathOperator{\cut}{cut}
\DeclareMathOperator{\disc}{disc}
\DeclareMathOperator{\Ig}{Ig}
\DeclareMathOperator{\IG}{IG}
\DeclareMathOperator{\Id}{id}
\DeclareMathOperator{\Rad}{Rad}
\DeclareMathOperator{\Sp}{Sp}

\renewcommand{\epsilon}{\varepsilon}
\renewcommand{\phi}{\varphi}

\begin{document}

\title[Igusa-type Functions and functional Equations]{Igusa-type
  Functions associated to finite formed Spaces and their functional
  Equations}

\author{Benjamin Klopsch}
\address{Mathematisches Institut \\ Heinrich-Heine-Universit\"at \\ D-40225
  D\"usseldorf \\ Germany}
\curraddr{Department of Mathematics \\ Royal Holloway,
  University of London \\ Egham TW20 0EX \\ United Kingdom}
\email{Benjamin.Klopsch@rhul.ac.uk}

\author{Christopher Voll}
\address{Max-Planck-Institut f\"ur Mathematik \\ Vivatsgasse 7 \\ D-53111
  Bonn \\ Germany}
\curraddr{School of Mathematics \\ University of Southampton \\
  Southampton SO17 1BJ \\ United Kingdom}
\email{C.Voll.98@cantab.net}
\thanks{The results in this paper form part of the first author's
  Habilitation thesis at the University of D\"usseldorf.
  The second author acknowledges support by the Deutsche
  Forschungsgemeinschaft and the Max-Planck-Gesellschaft.  He
  gratefully acknowledges the hospitality of the
  Heinrich-Heine-Universit\"at in D\"usseldorf and the
  Max-Planck-Institut f\"ur Mathematik in Bonn during the writing of
  this paper. This paper forms part of his Habilitation thesis at the
  University of D\"usseldorf.  }



\date{August 7, 2006 and, in revised form, October 25, 2007}


\begin{abstract}
  We study symmetries enjoyed by the polynomials enumerating
  non-degenerate flags in finite vector spaces, equipped with a
  non-degenerate alternating bilinear, hermitian or quadratic form.
  To this end we introduce Igusa-type rational functions encoding
  these polynomials and prove that they satisfy certain functional
  equations.
  
  Some of our results are achieved by expressing the polynomials in
  question in terms of what we call parabolic length functions on
  Coxeter groups of type~$A$. While our treatment of the orthogonal
  case exploits combinatorial properties of integer compositions and
  their refinements, we formulate a precise conjecture how in this
  situation, too, the polynomials may be described in terms of
  parabolic length functions.
\end{abstract}

\maketitle


\section{Introduction and discussion of results}

\subsection{General introduction}
Counting functions occurring naturally in algebra and geometry
frequently display symmetries which manifest themselves in the form of
certain functional equations. A classical example is the functional
equation satisfied by the Weil zeta function of a smooth, projective
algebraic variety defined over a finite field. More recently, Stanley
established similar symmetries for the Hilbert-Poincar\'e series of
graded algebras, with remarkable applications to counting problems in
combinatorics and topology; see \cite{Stanley/78}. The symmetries of
the Weil zeta functions lie at the heart of Denef and Meuser's proof
of a functional equation for certain $p$-adic integrals, called
Igusa's local zeta functions; see \cite{DenefMeuser/91}. The
phenomenon of functional equations also arises in the context of zeta
functions associated to groups and rings. Indeed, a recent result of
the second author establishes functional equations for local zeta
functions of finitely generated nilpotent groups, enumerating the
numbers of prime-power index subgroups; see \cite{Voll/06a}. It
depends on generalisations of methods developed for studying Igusa's
local zeta functions and draws on some of the ideas introduced in the
current paper.\footnote{The last two sentences were added during the
revision of the manuscript in October 2007. In the sequel we inserted
references to~\cite{Voll/06a} and \cite{duSG/06} where appropriate.}

In the current paper we are concerned with a counting problem in
geometric algebra: we study the numbers of flags of
\emph{non-degenerate} subspaces in a finite vector space, equipped
with an alternating bilinear, hermitian or quadratic form. The
polynomials giving these numbers can easily be computed, e.g.\
following Artin's classical book \cite{Artin/57} on Geometric
Algebra. Combining ideas from the theory of zeta functions, Coxeter
groups and combinatorics, we are able to establish remarkable
symmetries satisfied by these numbers which are far from evident from
the formulae. This is achieved by proving functional equations for
rational functions encoding the polynomials in question; see
Theorems~\ref{theorem_A} and \ref{theorem_B} below. A recurrent idea
in the present paper is to describe these polynomials in terms of
Coxeter groups and then to deduce functional equations from
generalisations of arguments of Igusa. In an important special case of
Theorem~A we can only conjecture such a description, leading to
Conjecture~\ref{conjecture_C}, which is of independent interest. The
proof we give for this special case of Theorem~\ref{theorem_A} is
combinatorial.

Our results have a precedent in the work of Igusa. In \cite[Part
II]{Igusa/89} Igusa establishes a formalism for studying certain
$p$-adic integrals associated to reductive algebraic groups. Igusa
computes a closed formula for these integrals in terms of the
associated Weyl groups and their root systems. In the case of
classical groups the data involved encode the number of flags of
\emph{totally isotropic} subspaces in finite polar spaces. In this
sense our work is complementary to that of Igusa. In both cases, the
functional equations may be understood in terms of the symmetry which
is given by (right-)multiplication by the longest element in a Coxeter
group.

Igusa's functions are closely related to $p$-adic integrals associated
to zeta functions of groups and rings, where many instances of
functional equations similar to the ones described in
Theorem~\ref{theorem_A} occur. The analytic properties of Euler
products of $p$-adic integrals of this type are also objects of
intense study. We refer the reader to~\cite{duS-ennui/03, duSG/00,
duSLubotzky/96, Voll/05, duSG/06, Voll/06a} for more information on
analytic properties of zeta functions of groups and their functional
equations.  Whilst we first encountered some of the Igusa-type
rational functions studied in the current paper in the context of zeta
functions of groups, they themselves are not generating functions.
Rather than encoding infinite arithmetic sequences they are associated
to finite formed spaces. We remark that a priori the rational
functions studied in the current paper do not have a natural
interpretation as $p$-adic integrals.

\subsection{Detailed statement of results}
\subsubsection{Theorem~\ref{theorem_A}}
In order to give a detailed statement of our results, we require some
notation. We fix a natural number~$n \in \N$ and consider
$n$-dimensional vector spaces~$V$ over a finite field~$F$, equipped
with a non-degenerate
\begin{itemize}
\item alternating bilinear form~$B$ (the `symplectic
  case'), 
\item hermitian form $B$ (the `unitary case') or 
\item quadratic
  form $f$ (the `orthogonal case').
\end{itemize}
In the symplectic and unitary cases, we formally define $f : V
\rightarrow F$ by $f(x) := B(x,x)$.  In the orthogonal case, we let
$B$ denote the bilinear form obtained by polarising $f$: if $\cha F
\ne 2$, then $B$ is non-degenerate symmetric, whereas, if $\cha F =
2$, then $B$ is alternating and possibly degenerate. The triple
$\calV := (V,B,f)$ will be called a \emph{formed space}. We also
introduce a parameter $\e$ equal to $1$ in the unitary case and equal
to $1/2$ otherwise; with this convention $F = \mathbb{F}_{q^{2 \e}}$
for a prime power $q$.

Recall that, by the classification of finite formed spaces (cf.,
e.g.,~\cite[Section~3.3]{Cameron/91}), $\calV$ decomposes as an
orthogonal direct sum of a certain number of hyperbolic planes and an
anisotropic space of dimension~$d \in \{0,1,2\}$. In the orthogonal
case we attach a sign~$\epsilon \in \{-1,1\}$ to $\calV$ if $n$ is
even, according to whether $d$ equals~$0$ or~$2$. The six
possibilities are given by the following table.

\begin{figure}[H]
\begin{center}
\begin{tabular}{|c|c!{\vrule width 1pt}c|c|c|c|}\hline
geometric type & $n$ & $d$ & $\epsilon$ & $\e$ \\\hline
symplectic&$2m$&$0$&--&$1/2$\\
unitary&$2m$&$0$&--&$1$\\
unitary&$2m+1$&$1$&--&$1$\\
orthogonal&$2m$&$0$&$1$&$1/2$\\
orthogonal&$2m+1$&$1$&--&$1/2$\\
orthogonal&$2m$&$2$&$-1$&$1/2$\\
\hline
\end{tabular}
\end{center}
\end{figure}
 
In the current paper we study rational functions incorporating the
numbers of $F$-rational points of the varieties of flags of
non-degenerate subspaces in~$\calV$. We write~$[n-1]$
for~$\{1,\dots,n-1\}$. By a \emph{non-degenerate flag of type $J =
  \{j_1,\dots,j_s\} \subseteq [n-1]$}, where $j_1<j_2<\dots<j_s$, we
mean a family $\bfU_J=(U_j)_{j\in J}$ of non-degenerate subspaces
of~$\calV$ with $U_{j_1}\subset\dots\subset U_{j_s}$ and $\dim U_j =j$
for each~$j\in J$. Let
$$
a_{\calV}^J(q):=|\{\bfU_J \mid \bfU_J\text{ non-degenerate flag of
  type }J\}|.
$$ Then $a^J_{\calV}(q)$ is a monic polynomial in $q$ (cf.~the remarks
at the end of this subsection regarding the orthogonal case), and we
set
$$
\alpha^J_{\calV}(q^{-1}):=a^J_{\calV}(q)/q^{\deg_qa^J_{\calV}}.
$$

We encode these numbers in rational functions as follows. Let $\bfX =
(X_i)_i$ be a finite family of independent indeterminates. Fix a
family of rational functions $\bfF = \left( F_{J}(\bfX) \right)_{J
  \subseteq [n-1]}$ in~$\bfX$ with the \emph{inversion property} that
\begin{equation}\tag{IP}\label{eq_IP}
  \text{for all }~I\subseteq[n-1]: \; F_{I}(\bfX^{-1}) = (-1)^{|I|}
  \sum_{J\subseteq I} F_J(\bfX). 
\end{equation}
A simple and naturally occurring
example~(cf.~\cite[Part~II]{Igusa/89}) of a family with this property
is $\left(\prod_{j\in J}\frac{X_j}{1-X_j}\right)_{J\subseteq [n-1]}.$
By defining
$$
\Ig_{\calV}(q^{-1},\bfX) := \Ig_{\calV,\bf{F}}(q^{-1},\bfX) :=
\sum_{J\subseteq[n-1]}\alpha^J_{\calV}(q^{-1})F_J(\bfX)
$$
we associate to $\calV$ and $\bfF$ a rational function in $q$ and the
variables $(X_i)_i$. Some explicit examples of these \emph{Igusa-type
  functions} may be found in the Appendix.

The first main result of this paper is

\begin{thmx} \label{theorem_A} 
  For each $n$-dimensional, non-degenerate formed space~$\calV$ the
  associated Igusa-type function satisfies the functional equation
  $$\Ig_{\calV}(q,\bfX^{-1})=(-1)^aq^{b}\,\Ig_{\calV}(q^{-1},\bfX),$$
  where the integers $a$ and $b$ are given by the table below
  \textup{(}with $m:=\lfloor\frac{n}{2}\rfloor$\textup{)}.

  \begin{figure}[H]
  \begin{center}
  \begin{tabular}{|c|c|c!{\vrule width 1pt}c|c|}\hline
  \textup{geometric type} & $n$ & $\epsilon$ & $a$ & $b$ \\\hline
  \textup{symplectic} &$2m$&--&$m-1$&$m(m-1)$\\
  \textup{unitary} &$n$&--&$\binom{n}{2}+n-1$&$\binom{n}{2}$\\
  \textup{orthogonal} &$2m$&$1$&$m+1$&$m^2$\\
  \textup{orthogonal} &$2m+1$&--&$m$&$m(m+1)$\\
  \textup{orthogonal} &$2m$&$-1$&$m$&$m^2$\\
  \hline
  \end{tabular}
  \end{center}
  \end{figure}
\end{thmx}

As we explained in the general introduction, we see
Theorem~\ref{theorem_A} primarily as a result about the polynomials
$\alpha^J_\calV(q^{-1})$; the choice of the family of rational
functions $\bf{F}$ is secondary. The inversion property \eqref{eq_IP}
satisfied by the rational functions $\bf{F}$ is a key ingredient which
we require for our subsequent combinatorial and group theoretical
considerations. We remark that rational functions $\bf{F}$ satisfying
the inversion property arise naturally in Igusa's work as well as in
the context of zeta functions of groups and rings; cf.~\cite{Voll/06,
Voll/06a}.

\begin{remark}[Analogy with Igusa's work]
  The analogy with Igusa's paper~\cite{Igusa/89} is the
  following. On~\cite[p.~706]{Igusa/89} Igusa gives a formula for a
  $p$-adic integral, essentially of the form
  $$
  Z(s)=\frac{\sum_{w\in W} q^{-l(w)}\prod_{\alpha_j \in w(R^-)} q^{A_j
      - B_js}}{\prod_{j=1}^\ell (1 - q^{A_j-B_js})},
  $$ where $q$ is a prime power, $W$ is a Weyl group, $l$ denotes the
  standard Coxeter length function, $S
  =\{\alpha_1,\dots,\alpha_\ell\}$ constitutes a basis for the root
  system, $R^-$ denotes the set of negative roots with respect to $S$,
  the parameters $A_j,B_j$ are suitable integers and $s$ is a complex
  variable. It is immediate that
  \begin{equation*}
    Z(s)=\sum_{J\subseteq[\ell]} \beta_J(q^{-1}) F_J(\bfX),
  \end{equation*}
  where $$\beta_J(q^{-1}) := \sum_{\substack{w\in W\\D(w)\subseteq I}}
  q^{-l(w)},\quad D(w) := \{j \in [\ell] \mid \alpha_j \in w(R^-)\},$$
  and $$F_J(\bfX) := \prod_{j\in J}\frac{X_j}{1-X_j},\quad X_j :=
  q^{A_j-B_js}.$$ The sets $D(w)$ may be interpreted as descent sets
  (cf.~Section~\ref{section_coxeter}). If $W$ comes from a classical
  group, the polynomials $\beta_J(q^{-1})$ which arise in this way
  carry a geometric meaning: the number $b_J(q)$ of flags of
  \emph{totally isotropic} subspaces in an associated finite polar
  space equals $\beta_J(q)$. In this case $b_J(q)$ gives the number of
  $\mathbb{F}_q$-points of a smooth projective variety, and we have
  $b_J(q) / q^{\deg_q(b_J)} = b_J(q^{-1}) = \beta_J(q^{-1})$.  In this
  sense the rational functions $\text{Ig}_\calV$ studied in the
  current paper, which are built in a similar way from the numbers of
  \emph{non-degenerate} flags, complement Igusa's function~$Z(s)$.

  The key to Igusa's functional equation is to interpret the inversion
  of the `variable' $q$ in terms of a natural symmetry of the root
  system of the Weyl group $W$. This symmetry arises from
  (right-)multiplication by the longest element $w_0 \in W$. A
  recurrent theme of the present paper is to give a suitable
  description of the polynomials $\alpha^J_\calV(q^{-1})$ in terms of
  Coxeter groups. Based on such a description one can then follow
  Igusa's approach to derive functional equations. We note that,
  contrary to the situation studied by Igusa, the `normalisation' of
  the polynomials $a^J_\calV(q)$ is indispensable -- unlike their
  counterparts $\beta_J$ and $b_J$, the polynomials $\alpha^J_\calV$
  and $a^J_\calV$ are typically not equal; cf.\ the examples given in
  the Appendix. Moreover, it is worth noting that as a side effect of
  passing to the normalised polynomials $\alpha^J_{\calV}(q^{-1})$ the
  assumption that $\calV$ is non-degenerate means no loss of
  generality.
\end{remark}
 
We now discuss the proof of Theorem~\ref{theorem_A}. In the symplectic
and unitary case it follows from Witt's Extension Theorem that the
respective isometry group acts transitively on the non-degenerate
flags of a given type. A simple stabiliser computation reveals that
the polynomials $\alpha^J_{\calV}(q^{-1})$ may be expressed in terms
of Gaussian polynomials (or $q$-binomial coefficients), which in turn
admit a well-known description in terms of the \emph{length function}
on a Coxeter group of type~$A_{\e n-1}$. The functional equation then
follows with the same argument which Igusa has given
in~\cite[Part~II]{Igusa/89}. It rests on the fact that, in a Coxeter
group, the effect of right-multiplication by the longest element on an
element's length and descent set is well understood.
 
In the orthogonal case, however, things are more intricate. To begin
with, the $a^J_{\calV}(q)$ flags of type $J$ come in up to $2^{|J|}$
isomorphism types and counting them together seems to be crucial for
the occurrence of a functional equation. But of course the natural
action of the respective orthogonal group on these flags is not
transitive.  The proof we give for this case of
Theorem~\ref{theorem_A} is based on a combinatorial analysis of the
polynomials~$\alpha^J_{\calV}(q^{-1})$ in terms of \emph{integer
  compositions} and their refinements. Complementing this approach, we
propose in Conjecture~\ref{conjecture_C} an explicit formula which
expresses these polynomials, too, in terms of Coxeter group data.

\begin{remark}[The orthogonal case in characteristic $2$] As is
  well-known, quadratic forms are intimately related to symmetric
  bilinear forms. In fact, over a field of characteristic not equal to
  $2$, the two notions lead to one and the same theory: a quadratic
  space $\calV = (V,B,f)$ over a field $F$ with $\cha F \ne 2$ can
  equally well be regarded as a symmetric bilinear space and vice
  versa. Such a space $\calV$ is said to be \emph{non-degenerate} if
  the bilinear form $B$ is non-degenerate, i.e.\ if the radical
  $\Rad(B) := \{ x \in V \mid \forall y \in V: B(x,y) = 0 \}$ is the
  zero subspace. In particular, enumerating non-degenerate flags in a
  quadratic space $\calV$ is the same as counting non-degenerate flags
  in the symmetric bilinear space $\calV$.
  
  In characteristic $2$, however, one has to distinguish more
  carefully between quadratic and symmetric bilinear forms. It is
  noteworthy that the analogous statement of Theorem~\ref{theorem_A}
  for symmetric bilinear spaces does not hold in characteristic $2$:
  in the Appendix we display a $4$-dimensional non-degenerate
  symmetric bilinear space whose associated `Igusa-type' function does
  not satisfy a functional equation.
  
  Now consider quadratic spaces $\calV = (V,B,f)$ over a field $F$
  with $\cha F = 2$. In this context $B$ is alternating and carries
  less information than $f$. There are basically two notions of
  `non-degeneracy', but unfortunately no standard terminology;
  cf.~\cite{Cameron/91}, \cite{Grove/02},
  \cite[Appendix~1]{MilnorHusemoller/73}, \cite{Pfister/95}. In this
  paper, we call $\calV$ \emph{non-defective} if the associated
  bilinear form $B$ is non-degenerate, i.e.~if the radical $\Rad(B) :=
  \{ x \in V \mid \forall y \in V: B(x,y) = 0 \}$ is the zero
  subspace. This can be thought of as a strong version of
  `non-degeneracy'; in particular, every non-defective quadratic space
  is even-dimensional. But enumerating non-defective flags in a
  quadratic space over $F$ is the same as counting non-degenerate
  flags in the induced alternating bilinear space, so we gain nothing
  new.  We call a quadratic space $\calV$ \emph{non-degenerate} if the
  restriction of $f$ to the radical $\Rad(B)$ is anisotropic, i.e.~if
  for all $x \in \Rad(B)$ either $x = 0$ or $f(x) \ne 0$. This concept
  of `non-degeneracy' is more flexible; in particular, there are
  non-degenerate quadratic spaces of any given dimension. Moreover,
  this turns out to be the right notion to formulate
  Theorem~\ref{theorem_A}. In fact, the polynomials $a_\calV^J(q)$
  counting non-degenerate flags of type $J$ are the same in all
  characteristics; see Section~\ref{section_orthogonal}.
\end{remark}

\subsubsection{Theorem~\ref{theorem_B}}
 
In the symplectic and unitary case, we prove a result which is
slightly more general than Theorem~\ref{theorem_A}. Rather than
counting flags which are non-degenerate with respect to a single
non-degenerate sesquilinear form $B$, we study the numbers of flags
which are non-degenerate with respect to a `flag of forms'. Loosely
speaking, a \emph{flag of sesquilinear forms~$\bfB$} of
type~$I\subseteq[n-1]$ is a family of sesquilinear forms such that
\begin{itemize}
\item all but the first form are degenerate,
\item each but the last form is defined on the radical of its
  successor, 
\item the last form is defined on the total space $V$ and 
\item the non-zero radicals constitute a flag of type $I$ in
  $V$.  
\end{itemize}
`Non-degeneracy' is defined inductively; see
Section~\ref{section_symplectic_unitary} for details.

Now let $\bfB$ be a sesquilinear flag of forms of type $I \subseteq
[n-1]$ on $V$. We denote by $a_{(V,\bfB)}^{J}(q)$ the number of flags
of type~$J \subseteq [n-1]$ which are non-degenerate with respect to
$\bfB$. In the symplectic case, both the type~$I$ of $\bfB$ and all
the sets~$J\subseteq[n-1]$ for which $a^J_{(V,\bfB)}(q)$ is non-zero
necessarily consist of even numbers. From the \emph{normalised}
polynomials
$$
\alpha_{(V,\bfB)}^J(q^{-1}) :=
a_{(V,\bfB)}^J(q)/q^{\deg_qa_{(V,\bfB)}^J}
$$
and a family of rational functions~${\bf
  F}=(F_{J}(\bfX))_{J\subseteq[n-1]}$ with the inversion
property~\eqref{eq_IP} we define, similarly as
above, a rational function
$$
\Ig_{(V,\bfB)}(q^{-1},\bfX) := \Ig_{(V,\bfB),{\bf
    F}}(q^{-1},\bfX)=\sum_{J\subseteq[n-1]}\alpha_{(V,\bfB)}^J(q^{-1})F_J(\bfX).
$$
The second main result of this paper is

\begin{thmx}\label{theorem_B}
  For each $n$-dimensional vector space~$V$, equipped with a flag of
  alternating bilinear \textup{(}respectively hermitian\textup{)}
  forms~$\bfB$ of type $I=\{i_1,\dots,i_r\}_<\subseteq[n-1]$, the
  associated Igusa-type function satisfies the functional equation
  \begin{equation*}
   \Ig_{(V,\bfB)}(q,\bfX^{-1})=(-1)^a
   q^b\,\Ig_{(V,\widetilde{\bfB})}(q^{-1},\bfX),
  \end{equation*}
  where $\widetilde{\bfB}$ is a flag of forms of
  type~$\widetilde{I}:=\{n-i \mid i\in I\}$ and the integers $a$ and
  $b$ are given by the table below \textup{(}with $m :=
  \lfloor\frac{n}{2}\rfloor$\textup{)}.

  \begin{figure}[H]
  \begin{center}
  \begin{tabular}{|c!{\vrule width 1pt}c|c|}\hline
   \textup{geometric type} & $a$ & $b$ \\ \hline
   \textup{symplectic} &$m-1$&$m(m-1)+((i_2-i_1)i_1+\dots+(n-i_r)i_r)/2$\\
   \textup{unitary} &$n-1+b$&$\binom{n}{2} + (i_2-i_1)i_1+\dots+(n-i_r)i_r$\\
   \hline\end{tabular}
  \end{center}
  \end{figure}
\end{thmx}

Note that, for $I=\varnothing$, Theorem~\ref{theorem_B} specialises to
Theorem~\ref{theorem_A} in the symplectic and unitary case,
respectively.

To prove Theorem~\ref{theorem_B} we show that the
functions~$\alpha_{(V,\bfB)}^J(q^{-1})$ are polynomials which may be
described in terms of a certain \emph{statistic} on the Coxeter group
$W$ of type~$A_{\e n -1}$. This statistic associates to an element $w
\in W$ the sum of its ordinary length~$l(w)$ with respect to the
standard Coxeter generating set~$S = \{s_1,\dots,s_{\e n-1}\}$ and its
`parabolic length'~$l_\upL^{(\e \widetilde{I})^c}(w)$. The
\emph{parabolic length} $l_\upL^{(\e \widetilde{I})^c}(w)$ is the
Coxeter length of the distinguished representative of shortest length
in the left coset~$w W_{(\e \widetilde{I})^c}$ of the standard
parabolic subgroup $W_{(\e \widetilde{I})^c} = \langle s_i\in S \mid
\e n - i \not \in \e I \rangle$.
 
In fact, in Section~\ref{section_symplectic_unitary} we show that
Theorem~\ref{theorem_B} can be deduced from Theorem~\ref{theorem_1}, a
general result on rational functions defined in terms of linear
combinations of parabolic length functions and characters on certain
subgroups of finite Coxeter groups. Indeed, Theorem~\ref{theorem_1}
extends to a slightly more general setting Igusa's key idea to deduce
functional equations from features of the map induced by (right-)
multiplication by the longest element.

Our initial interest in the Igusa-type functions $\Ig_{(V,\bfB)}$
arose from our study of the zeta functions counting subgroups of
higher Heisenberg groups.  In~\cite{KlopschVoll/05} we introduce an
equivalence relation, coarser than homothety, on the set of complete
$\mathbb{Z}_p$-lattices in a non-degenerate symplectic $p$-adic vector
space~$\mathbb{Q}_p^{2m}$ such that equivalence classes of lattices
are in one-to-one correspondence with the vertices of the affine
Bruhat-Tits building for the symplectic
group~$\text{Sp}_{2m}(\mathbb{Q}_p)$. For a flag of forms $\bfB$ of
type $I=\{2i\}$, $i\in[m]$, and suitable choices of $\bfF$ the Igusa
functions~$\Ig_{(\F_q^{2m},\bfB),\bfF}$ may be regarded as generating
functions, enumerating lattices in an equivalence class indexed by a
special vertex of type~$i$. We refer to~\cite{KlopschVoll/05} for
details.

\subsubsection{Conjecture~\ref{conjecture_C}}

In the last part of the paper we formulate a precise conjecture
describing the polynomials $\alpha^J_{\calV}(q^{-1})$ in the
orthogonal case.  If it holds, the orthogonal case of
Theorem~\ref{theorem_A} also follows from Theorem~\ref{theorem_1}.
Moreover, a proof of Conjecture~\ref{conjecture_C} would constitute a
first step towards extending Theorem~\ref{theorem_B} to the orthogonal
case.

We introduce the subgroup~$\mathcal{C}_n$ of `chessboard elements' in
the symmetric group $\sym_n$ on~$n$ letters. A permutation is a
\emph{chessboard element} if the non-zero entries of its associated
permutation matrix all fit either on the black or on the white squares
of an $n\times n$-chessboard. In Section~\ref{section_conjecture} we
define linear characters~$\chi_{\epsilon}$ on~$\mathcal{C}_n$ and a
certain linear combination~$L$ of parabolic length functions
on~$\sym_n$. By $D_\upL(w)$ we denote the left-descent set of the
permutation~$w$ (cf.\ Section~\ref{section_coxeter}).

\begin{conjecture}\label{conjecture_C} For each 
  $n$-dimensional, non-degenerate quadratic space $\calV$ and each $J
  \subseteq [n-1]$,
  $$
  \alpha^J_{\calV}(q^{-1}) = \sum_{\substack{w \in \mathcal{C}_n \\ 
      D_\upL(w) \subseteq J}} \chi_{\epsilon}(w) q^{-L(w)}.
  $$
\end{conjecture}

\subsection{Organisation and Notation}

The structure of the paper is as follows. In
Section~\ref{section_coxeter} we derive functional equations for
rational functions defined in terms of parabolic length functions on
Coxeter groups.  Theorem~\ref{theorem_1}, the main result of
Section~\ref{section_coxeter}, is applied to prove
Theorem~\ref{theorem_B} in Section~\ref{section_symplectic_unitary}.
In Section~\ref{section_orthogonal} we prove the orthogonal case of
Theorem~\ref{theorem_A}.  In Section~\ref{section_conjecture} we give
a more precise statement of Conjecture~\ref{conjecture_C}. Some
explicit examples of Igusa-type functions can be found in the
Appendix.

\bigskip

\noindent We use the following notation.
\nopagebreak

\medskip

\begin{tabular}{l|l}
  $\N$ & the set $\{1,2,\dots\}$ of natural numbers \\
  $I_0$ & the set $I\cup\{0\}$ for $I\subseteq\N$ \\
  $[a,b]$ & the interval $\{a,a+1,\dots,b\}$ for integers $a,b$ \\
  $[a]$ & the interval $[1,a]$ for an integer $a$ \\
  $\{i_1,\dots,i_r\}_{<}$ & the set $\{i_1,\dots,i_r\} \subseteq \N_0$ with
  $i_1<\dots<i_r$ \\
  $x I$ & the set $\{ x i  \mid  i \in I \}$ for $ I\subseteq
  \N$ and a rational number $x$ \\
  $I^c$ &  the set $[n-1]\setminus I$ for $I\subseteq[n-1]$, \\
  & \quad where $n$ is clear from the context\\
  $\widetilde{I}$ & the set $\{n-i \mid  i\in I\}$ for $I\subseteq[n-1]$, \\
  & \quad where $n$ is clear from the context\\
  $I-\mathbf{t}$ & the set $\{i_1-t_1,\dots,i_r-t_r\} \cap \N$ for
  $I=\{i_1,\dots,i_r\}_{<} \subseteq \N$\\  
  & \quad and $\mathbf{t} = (t_1,\dots,t_r) \in \N_0^r$ \\
  $I-j$ & the set $\{i-j \mid i \in I \} \cap \N$ for $I \subseteq \N$,
  $j \in \N_0$\\
  $\binom{a}{b}$ & the ordinary binomial coefficient for $a, b \in \N_0$\\
  $\binom{a}{b}_{\! X}$ & the polynomial
  $\prod_{i = 0}^{b-1} (1 - X^{a-i})/(1 - X^{b-i})$, \\
  & \quad where $a,b \in \N_0$ with $a \geq b$ \\
  & Note: The \emph{$q$-binomial coefficient} $\binom{a}{b}_{\! q}$ gives \\
  & \quad the number of subspaces of dimension $b$  in $\F_q^a$.\\
  $\binom{n}{J}_{\! X}$ & the polynomial
  $\binom{n}{j_s}_{\! X} \binom{j_s}{j_{s-1}}_{\! X} \dots
  \binom{j_2}{j_1}_{\! X}$,\\
  & \quad where $J = \{j_1,\dots,j_s\}_< \subseteq[n-1]_0$ for $n \in \N$\\
  & Note: $\binom{n}{J}_{\! q} \text{ gives the number of flags of
    type~$J\setminus\{0\}$ in~$\F_q^n$.}$\\
  $\lfloor x \rfloor$ & the greatest integer not exceeding the
  rational number $x$\\
  $\mathcal{P}(S)$ & the power set of a set $S$\\
  $\sym_n$ & the symmetric group on $n$ letters.
\end{tabular}

\bigskip 

\noindent Throughout this paper $n \in \N$ and $m = \lfloor n/2 \rfloor$. We
shall write
\begin{align*}
  I & = \{i_1,\dots,i_r\}_<,& J & = \{j_1,\dots,j_s\}_< & & \text{ for
    subsets of $[n-1]$ or $[n]$, and} \\
  G & = \{g_1,\dots,g_k\}_<,& H & & & \text{ for
    subsets of $[m]$.}
\end{align*}


\section{Rational functions from Coxeter groups}\label{section_coxeter}

In this section we prove functional equations for a family of rational
functions associated to finite Coxeter systems.
Theorem~\ref{theorem_B} will turn out to be a consequence of
Theorem~\ref{theorem_1}, the main result of the current Section.

Let $(W,S)$ be a finite Coxeter system of rank $n-1$ with root
system~$\Delta$. To ease notation we will frequently identify the set
of Coxeter generators $S=\{s_1,\dots,s_{n-1}\}$ with the set of
integers~$[n-1]$. For each~$I\subseteq S$ we denote by $W_I$ the
corresponding standard parabolic subgroup of~$W$ generated by the
elements in~$I$ and by~$\Delta_I$ the induced root system. We denote
by $l$ the length function on $W$ with respect to~$S$. The length of
an element~$w$ may either be interpreted as the length of a shortest
word in the elements of~$S$ representing the group element or as the
number of positive roots that are sent to negative roots by~$w$. The
group~$W$ has a unique longest element~$w_0$, whose length equals
$|\Delta|/2$. It is well-known (cf.~\cite[Section~1.8]{Humphreys/90})
that, for each~$w\in W$,
$$
l(w_0w)+l(w) = l(ww_0)+l(w) = l(w_0).
$$

The rational functions studied in this section are defined in terms of
more general length functions. For each $I\subseteq S$ set
\begin{align*}
  W^I_{\upL}& := \{w\in W \mid \forall s\in I:\,l(ws)>l(w)\},\\
  W^I_{\upR}& := \{w\in W \mid \forall s\in I:\,l(sw)>l(w)\}.
\end{align*}
We will need the following lemma
(\cite[Proposition~2.1.7]{Scharlau/95}).

\begin{lemma}\label{lemma_1}
  Let $I \subseteq S$. Then $W^I_{\upL}$ \textup{(}respectively
  $W^I_{\upR}$\textup{)} is a left \textup{(}respectively
  right\textup{)} transversal to $W_I$ in $W$, i.e.\ for every $w \in
  W$ there are unique elements
  $$
  u_\upL\in W^I_{\upL}, \; v_\upL\in W_I \; \text{ and } \; u_\upR\in
  W^I_{\upR}, \; v_\upR\in W_I
  $$
  such that
  $$
  w = u_\upL v_\upL = v_\upR u_\upR.
  $$
  In particular, $u_\upL$ is the unique element of shortest length in
  the left coset $w W_I$ and $u_\upR$ is the unique element of shortest
  length in the right coset $W_I w$. Moreover,
  $$
  l(w) = l(u_\upL)+l(v_\upL)=l(v_\upR)+l(u_\upR).
  $$
  The elements $u_\upL \in w W_I$ and $u_\upR \in W_I w$ are also
  characterised by the fact that they send positive roots
  of~$\Delta_I$ to positive roots.
\end{lemma}

\begin{definition}[Parabolic length]
  For each $I\subseteq S$ and $w\in W$ we set
  \begin{align*}
    l^I_\upL(w)&:= l(u_\upL),\\
    l^I_\textup{R}(w)&:= l(u_\upR).
  \end{align*}
  We call $l^I_\upL$ (respectively $l^I_\textup{R}$) the \emph{left}
  (respectively \emph{right}) \emph{parabolic length function} on~$W$
  associated to~$I$. We write $\bfl_L:=(l^I_\upL)_{I\subseteq S}$ and
  $\bfl_\upR:=(l^I_\textup{R})_{I\subseteq S}$.
\end{definition}
Note that for $I = \varnothing$ the corresponding parabolic length
functions reduce to the ordinary Coxeter length function:
$l^{\varnothing}_\upL=l^{\varnothing}_\upR=l$. Moreover,
$l^{S}_\upL=l^{S}_\upR=0$.

\begin{lemma} \label{lemma_2}
  For each $I\subseteq S$ and $w\in W$ we have
  \begin{align}
    l^I_\upL(w_0w)+l^I_\upL(w) = l^I_\upL(w_0), & \quad
    l^I_\upL(ww_0)+l_\upL^{I^{w_0}}(w) = l^I_\upL(w_0),\label{eq_1}\\
    l^I_\upR(w_0w)+l_\upR^{I^{w_0}}(w) = l^I_\upR(w_0), & \quad
    l^I_\upR(ww_0)+l^I_\upR(w) = l^I_\upR(w_0).\label{eq_2}
 \end{align}
\end{lemma}

\begin{proof}
  Let $v_0$ denote the longest element in~$W_I$.  Then
  $$
  l^I_\upL(w_0)=l(w_0)-l(v_0)=l(w_0)- |\Delta_I|/2.
  $$
  Write $w=u_\upL v_\upL$ as in Lemma~\ref{lemma_1}. We may then write
  $$
  w_0=u'v'v_\upL^{-1}u_\upL^{-1}
  $$
  with $u'\in W$, $v'\in W_I$ such that $v'v_\upL^{-1}=v_0$. It
  follows that
  \begin{align*}
    l(u') &= l(w_0) -
    l(v'v_\upL^{-1}u_\upL^{-1}) = l(w_0)-l(u_\upL v_0)\\
    & = l(w_0) - l(u_\upL) - l(v_0) = (l(w_0) - l(v_0)) -
    l(u_\upL) \\
    & = l^I_\upL(w_0) - l^I_\upL(w).
  \end{align*}
  Clearly, $w_0 w W_I = u' W_I$. But $u'$ sends positive roots of
  $\Delta_I$ to positive roots and is thus the unique coset
  representative of shortest length. Hence $l^I_\upL(w_0w) = l(u') =
  l^I_\upL(w_0) - l^I_\upL(w)$. This gives the first equation in
  \eqref{eq_1}.
  
  For the second equation in \eqref{eq_1}, note that conjugation by~$w_0$
  yields $l^I_\upL(w_0)=l_{\upL}^{I^{w_0}}(w_0)$ and thus
  
  $$
  l^I_\upL(ww_0) = l_{\upL}^{I^{w_0}}(w_0w) = l_\upL^{I^{w_0}}(w_0)
  - l_\upL^{I^{w_0}}(w) = l^I_\upL(w_0)-l_\upL^{I^{w_0}}(w).
  $$
  
  We omit the analogous proofs for the equations~\eqref{eq_2}.
\end{proof}

Another important invariant of an element of a Coxeter group which we
shall need is its (left) \emph{descent set}~$D_\upL(w):=\{s\in S \mid
l(sw)<l(w)\}$. Note that
\begin{equation}\label{eq_3}
D_\upL(ww_0)=D_\upL(w)^c:=\{s\in S \mid s\not\in D_\upL(w)\}.
\end{equation}

Elements of Coxeter groups of type~$A$ can be regarded as permutation
matrices. It is noteworthy that both the descent sets and the values
of the various parabolic length functions are easily read off from the
associated matrices.

\begin{lemma} \label{lemma_3}
  Let $(F_{J}(\bfX))_{J\subseteq S}$ be a family of rational functions
  with the inversion property~\eqref{eq_IP}.
  Then, for all $I\subseteq S$,
  $$\sum_{ I\subseteq J \subseteq
    S}F_{J}(\bfX^{-1})=(-1)^{|S|}\sum_{I^c\subseteq J \subseteq S}
  F_{J}(\bfX)$$
\end{lemma}
\begin{proof} 
  This is an easy calculation. See~\cite[Lemma~7]{Voll/06}.
\end{proof}

We now fix a family of rational functions~$\bfF = (F_J(\bfX))_{J
  \subseteq S}$ with the inversion property~\eqref{eq_IP} and an
independent indeterminate $Y$. We choose a family~$\bfb =
(b_I)_{I\subseteq S}$ of integers and define the statistics
$\bfb\cdot\bfl_\upL$ and~$\bfb\cdot\bfl_\upR$ on~$W$ by setting, for
$w\in W$,
$$
\bfb\cdot\bfl_\upL(w):=\sum_{I\subseteq S}b_Il^I_\upL(w), \qquad
\bfb\cdot\bfl_\upR(w):=\sum_{I\subseteq S}b_Il^I_\upR(w).
$$
Similarly, we write $\bfb^{w_0}\cdot\bfl_\upL$ and
$\bfb^{w_0}\cdot\bfl_\upR$ to denote the statistics associating to~$w$
the elements $\sum_{I\subseteq S} b_I l_{\upL}^{I^{w_0}}(w)$ and
$\sum_{I\subseteq S} b_I l_{\upR}^{I^{w_0}}(w)$, respectively. Let
$W'\subseteq W$ be a subgroup with $w_0\in W'$, and
$\chi:W'\rightarrow \mathbb{C}^*$ a (linear) character of~$W'$.

\begin{definition}
  With the given data we define the following rational functions:
  \begin{align*}
    \IG_{\upL}^{W',\bfb,\chi,\bfF}(Y,\bfX)& := \sum_{w\in W'} \chi(w)
    Y^{\bfb\cdot\bfl_\upL(w)} \sum_{\substack{D_\upL(w) \subseteq J
        \subseteq S}}F_J(\bfX), \\
    \IG_{\upR}^{W',\bfb,\chi,\bfF}(Y,\bfX)& := \sum_{w\in W'} \chi(w)
    Y^{\bfb\cdot\bfl_\upR(w)} \sum_{\substack{D_\upL(w) \subseteq J
        \subseteq S}}F_J(\bfX).
  \end{align*}
\end{definition}  

The main result of the current section is
\begin{theorem}\label{theorem_1} 
  The following functional equations hold:
  \begin{align}
    \IG_{\upL}^{W',\bfb,\chi,\bfF}(Y^{-1},\bfX^{-1}) & = (-1)^{|S|}
    \chi(w_0) Y^{-\bfb\cdot\bfl_\upL(w_0)}
    \IG_{\upL}^{W',\bfb^{w_0},\chi,\bfF}(Y,\bfX), \label{eq_4} \\
    \IG_{\upR}^{W',\bfb,\chi,\bfF}(Y^{-1},\bfX^{-1}) & =(-1)^{|S|}
    \chi(w_0) Y^{-\bfb\cdot\bfl_\upR(w_0)}
    \IG_{\upR}^{W',\bfb,\chi,\bfF}(Y,\bfX). \label{eq_5}
  \end{align}
\end{theorem}

\begin{proof}
  The equations
  \begin{equation} \label{eq_6}
  \begin{split}
    \bfb^{w_0}\cdot\bfl_\upL(ww_0) + \bfb\cdot\bfl_\upL(w) & =
    \bfb\cdot\bfl_\upL(w_0), \\
    \bfb\cdot\bfl_\upR(ww_0) + \bfb\cdot\bfl_\upR(w) & =
    \bfb\cdot\bfl_\upR(w_0)
  \end{split}
  \end{equation}
  are immediate consequences of Lemma~\ref{lemma_2}. Therefore,
  by~\eqref{eq_6}, by Lemma~\ref{lemma_3} and by~\eqref{eq_3},
  \begin{align*}
    \IG_{\upL}^{W',\bfb,\chi,\bfF} & (Y^{-1},\bfX^{-1}) = \sum_{w\in W'}
    \chi(w) Y^{-\bfb\cdot\bfl_\upL(w)} \sum_{\substack{D_\upL(w) \subseteq
        J \subseteq S}} F_J(\bfX^{-1}) \\
    = (-1)^{|S|} & \chi(w_0)^{-1} Y^{-\bfb\cdot\bfl_\upL(w_0)}
    \sum_{w\in W'} \chi(ww_0) Y^{\bfb^{w_0}\cdot\bfl_\upL(ww_0)}
    \sum_{\substack{D_\upL(ww_0)\subseteq
        J\subseteq S}} F_J(\bfX)\\
    = (-1)^{|S|} & \chi(w_0) Y^{-\bfb\cdot\bfl_\upL(w_0)}
    \IG_{\upL}^{W',\bfb^{w_0},\chi,\bfF}(Y,\bfX).
  \end{align*}
  The equation~\eqref{eq_5} is proved analogously.
\end{proof}

In this paper we shall see instances of both types of functional
equations presented in Theorem~\ref{theorem_1}. In
Section~\ref{section_symplectic_unitary} we demonstrate that
Theorem~\ref{theorem_B} is a consequence of~\eqref{eq_4}.  Note that
in the special case $\bfF=\left(\prod_{j\in
    J}\frac{X_j}{1-X_j}\right)_{J\subseteq[n-1]}$, replacing~$\bfb$
by~$\bfb^{w_0}$ in~\eqref{eq_4} simply amounts to inverting the order
of the variables $X_1,\dots,X_{n-1}$. If Conjecture~\ref{conjecture_C}
holds, the orthogonal case of Theorem~\ref{theorem_A} follows
from~\eqref{eq_5}.


\section{The symplectic and unitary case}\label{section_symplectic_unitary}

In this section we study the polynomials enumerating flags which are
non-dege\-nerate with respect to a `flag of sesquilinear forms'. Our
aim is to proof Theorem~\ref{theorem_B}. Let $V$ be an $n$-dimensional
vector space over a field~$F$.
Let~$I=\{i_1,\dots,i_r\}_<\subseteq[n-1]$, and set $i_0 := 0$,
$i_{r+1} := n$.

\begin{definition}[Flag of forms]
  We say that~$V$ is equipped with a \emph{flag of alternating
    bilinear} (respectively \emph{hermitian}) \emph{forms}
  $\bfB=(B_{i_1},\dots,B_{i_{r+1}})$ of type~$I$ if there is a
  filtration of subspaces
  $$
  \{0\}=:R_{i_0}\subset R_{i_1}\subset\dots\subset R_{i_r}\subset
  R_{i_{r+1}}:=V
  $$
  such that
  \begin{enumerate}
  \item[(a)] for all $i\in I$, $\dim R_i=i$;
  \item[(b)] for all $ \rho\in [r+1]$, $B_{i_\rho}$ is an alternating
    bilinear (respectively hermitian) form $B_{i_\rho}: R_{i_\rho}
    \times R_{i_\rho} \rightarrow F$ with
    $$
    \Rad(B_{i_\rho}):=\{x\in R_{i_\rho} \mid \forall y\in
    R_{i_\rho}:\; B_{i_\rho}(x,y)=0\}=R_{i_{\rho-1}}.
    $$
  \end{enumerate}
  We call the sequence $\bfR=(R_{i_1},\dots,R_{i_r})$ the
  \emph{flag of radicals} associated to the flag of
  forms~$\bfB$.
\end{definition}

Note that, given a flag of sesquilinear forms~$\bfB$ of type~$I$ on
$V$ with flag of radicals~$\bfR$ and $\rho\in[r+1]$, we have a flag of
forms $(B_{i_1},\dots,B_{i_\rho})$ of type $\{i_1,\dots,i_{\rho-1}\}$
on $R_{i_\rho}$ with flag of radicals $(R_{i_1}, \dots,
R_{i_{\rho-1}})$ and a flag of forms~$(\overline{B_{i_{\rho+1}}},
\dots, \overline{B_{i_{r+1}}})$ of type~$\{i_\varrho - i_\rho \mid
\rho < \varrho \leq r \}$ on $V/R_{i_\rho}$ with flag of radicals
$(R_{i_{\rho+1}}/R_{i_\rho},\dots,R_{i_r}/R_{i_\rho})$.

\begin{definition}[Non-degeneracy] \label{definition_4}
  Given a flag of sesquilinear forms~~$\bfB$ of type~$I$ on~$V$ with
  flag of radicals~$\bfR$ as above, we say that a subspace~$U\subseteq
  V$ is \emph{non-degenerate} with respect to $\bfB$ if for each
  $\rho\in[r+1]$, 
  \begin{enumerate}
  \item[(a)] $U\cap R_{i_\rho}$ is non-degenerate with respect to
    $(B_{i_1},\dots,B_{i_\rho})$ and
  \item[(b)] $(U+R_{i_\rho})/R_{i_\rho}$ is \emph{non-degenerate} with
    respect to
    $(\overline{B_{i_{\rho+1}}},\dots,\overline{B_{i_{r+1}}})$.
  \end{enumerate}
  
  A flag $\mathbf{U}_J = (U_j)_{j\in J}$ of subspaces of~$V$ of type
  $J\subseteq[n-1]$, i.e.\ an ascending chain of subspaces with $\dim
  U_j =j$ for each~$j\in J$, is said to be \emph{non-degenerate} with
  respect to $\bfB$ if each of its constituents~$U_j$ is.
\end{definition}

These definitions are illustrated by the following simple example.

\begin{example} 
  Suppose that $n$ is even, $r=1$ and $i_1\in[n-1]$ is even. Then a
  flag of alternating bilinear forms $\bfB = (B_{{i_1}},B_n)$ consists
  of a (degenerate) alternating bilinear form~$B_n$ on~$V$ with
  ${i_1}$-dimensional radical $R_{{i_1}}$, which in turn supports a
  non-degenerate form~$B_{i_1}$.  A subspace $U\subseteq V$ is
  non-degenerate with respect to~$\bfB=(B_{i_1},B_n)$ if $U\cap
  R_{i_1}$ is non-degenerate with respect to~$B_{i_1}$ and
  $(U+R_{i_1})/R_{i_1}$ is non-degenerate with respect to
  $\overline{B_n}$.
\end{example}

We shall now assume that $V$ is an $n$-dimensional vector space over a
finite field $F$, equipped with a flag of sesquilinear forms~$\bfB$ of
type~$I$. As in the introduction we write $\e=1/2$ in the symplectic
case and $\e=1$ in the unitary case so that $F = \mathbb{F}_{q^{2\e}}$
for some prime power $q$. Let $J\subseteq[n-1]$ and
define~$a^J_{(V,\bfB)}(q)$ to be the number of flags of type~$J$ which
are non-degenerate with respect to $\bfB$. We set
$$
\alpha_{(V,\bfB)}^J(q^{-1}) :=
a^J_{(V,\bfB)}(q)/q^{\deg_qa_{(V,\bfB)}^J}
$$
and shall frequently write $a^J_{n,I}(q)$ for
$a^J_{(V,\bfB)}(q)$ and $\alpha^J_{n,I}(q^{-1})$ for
$\alpha^J_{(V,\bfB)}(q^{-1})$.  Recall that, in the symplectic
case, both the type~$I$ of a flag of forms and all the
sets~$J\subseteq[n-1]$ for which
$a^J_{(V,\bfB)}(q)$ is non-zero consist necessarily of even
numbers.

\begin{definition} \label{definition_5} 
  Given a family $\bfF=(F_J(\bfX))_{J\subseteq[n-1]}$ of rational
  functions with the inversion property~\eqref{eq_IP} we define
  $$
  \Ig_{(V,\bfB)}(q^{-1},\bfX):=\Ig_{(V,\bfB),\bfF}(q^{-1},\bfX) =
  \sum_{J\subseteq[n-1]}\alpha_{(V,\bfB)}^J(q^{-1})F_J(\bfX).
  $$
\end{definition}

Theorem~\ref{theorem_B} states that these Igusa-type functions satisfy a
functional equation. In the remainder of the current section we show
how this can be deduced from the first assertion of
Theorem~\ref{theorem_1}. Fix a family of rational
functions~$\bfF=(F_J(\bfX))_{J\subseteq [n-1]}$ with the inversion
property~\eqref{eq_IP}, and define
$$
\e \bfF := (F_{\e^{-1} J'}(\bfX))_{J' \subseteq [\e n - 1]}.
$$
Let $W':= W :=\sym_{\e n}$ be the full symmetric group on $\e n$
letters, let $\chi$ be the trivial character on $W'$ and set, for each
$J' \subseteq [\e n-1]$,
\begin{equation*} 
  b_{J'} := \delta(\text{`}J' = \varnothing\text{'}) +
  \delta(\text{`}J' = \{s_i \mid i \not \in
  \e \widetilde{I} \}\text{'}), 
\end{equation*}
where the Kronecker-delta $\delta(E) \in \{1,0\}$ reflects whether or
not equation $E$ holds.

In order to derive Theorem~\ref{theorem_B}
from \eqref{eq_4} it suffices to show that
$$
\Ig_{(V,\bfB),\bfF}(q^{-1},\bfX) = \IG_{\upL}^{W',\bfb,\chi,\e \bfF}
((-q)^{-1/\e},\bfX)
$$
for the given data $W'$, $\bfb$, $\chi$ and $\e \bfF$. Clearly, it is enough
to prove

\begin{proposition}
  Let $J\subseteq[n-1]$ such that $\e J\subseteq[\e n-1]$. Then
  \begin{equation} \label{eq_7}
    \alpha_{n,I}^J(q^{-1}) = a_{n,I}^J(q) / {q^{2\e\deg_q{\binom{
    n}{J}_{\! q}}}} = \sum_{\substack{w\in \sym_{\e n} \\ D_\upL(w) \subseteq
    \e J}} Y^{\lambda_{n,I}(w)},
  \end{equation}
  where $Y:=(-q)^{-1/\e}$ and $\lambda_{n,I}(w) := l(w) + l^{(\e
    \widetilde{I})^c}_{\upL}(w)$ for all $w \in W$.
\end{proposition}

\begin{proof}
  We will first prove \eqref{eq_7} in the case $I =
  \varnothing$. The proof consists of a simple index computation in
  the respective isometry group, i.e.\ in the symplectic group
  $\Sp_{n}(\F_q)$ or the unitary group
  $\textup{U}_n(\mathbb{F}_{q^2})$. The proof in the general case is
  then based on a recursive expression for the numbers~$a_{n,I}^J(q)$.
  
  So assume that $I=\varnothing$. Then $\bfB = (B)$ simply specifies a
  non-degenerate alternating bilinear (respectively hermitian) form on
  $V$. The respective isometry group $\Sp_{n}(\F_q)$ or
  $\textup{U}_n(\mathbb{F}_{q^2})$ acts transitively on the
  non-degenerate flags of type~$J$, so it suffices to compute the
  stabiliser of any one of them. We construct a `standard'
  non-degenerate flag ${\bf U}_J:=(U_{j_1},\dots,U_{j_s})$ of type~$J
  = \{ j_1, \dots, j_s \}_<$ in the following way. In the symplectic
  case, choose a symplectic basis $E =
  (e_1,f_1,\dots,e_{n/2},f_{n/2})$ for~$V$ (i.e.\
  $B(e_i,f_j)=\delta_{ij}$, $B(e_i,e_j) = B(f_i,f_j)=0$) and set
  $U_j:=\langle e_1,f_1,\dots,e_{j/2},f_{j/2} \rangle$ for $j \in J$.
  In the unitary case, choose a unitary basis $E = (e_1,\dots,e_n)$
  for~$V$ (i.e.\ $B(e_i,e_j)=\delta_{ij}$) and set $U_j:=\langle
  e_1,\dots,e_j\rangle$ for $j \in J$. It is not difficult to verify
  that an element of the respective isometry group of $(V,B)$
  stabilises $\bfU_J$ if and only if its matrix $M_n$ with respect to
  the basis $E$ is of block diagonal form
  $$
  M_n =
  \left(
  \begin{array}{cccc}
  M_{j_1} &            &       &           \\
          & M_{j_2-j_1}&       &           \\
          &            & \ddots&           \\
          &            &       & M_{n-j_s} \\
  \end{array}
  \right)
  $$
  with $M_{j_\sigma-j_{\sigma-1}}$ in the respective smaller
  isometry group for all $\sigma\in[s+1]$, where $j_0:=0, j_{s+1}:=n$.
  Thus
  $$
  a^J_{n,\varnothing}(q) = 
  \begin{cases}
    |\text{Sp}_{n}(\F_q)| / \prod_{\sigma\in[s+1]}
    |\text{Sp}_{j_\sigma-j_{\sigma-1}}(\F_q)| & \text{in the
      symplectic case,} \\
    |\text{U}_{n}(\mathbb{F}_{q^2})| / \prod_{\sigma\in[s+1]}
    |\text{U}_{j_\sigma-j_{\sigma-1}}(\mathbb{F}_{q^2})| & \text{in
      the unitary case.}
  \end{cases}
  $$
  Employing the well-known formulae (cf.~\cite[p.~147]{Artin/57},
  \cite[Theorems~3.12 and~11.28]{Grove/02})
  \begin{align*}
    |\Sp_{n}(\F_q)| & = q^{\binom{n+1}{2}}\prod_{i\in[n/2]} (1-q^{-2i}),\\
    |\textup{U}_n(\mathbb{F}_{q^2})| & = q^{n^2}
    \prod_{i\in[n]}(1-(-q^{-1})^i)
  \end{align*}
  and using the notation $Y = (-q)^{- 1/ \e}$ we obtain
  $$
  \alpha_{n,\varnothing}^J(q^{-1}) = \frac{\prod_{i\in[\e n]} (1 -
    Y^i)}{\prod_{\sigma\in[s+1]} \prod_{\iota \in [ \e
      (j_\sigma-j_{\sigma-1})]} (1-Y^\iota)} = \binom{\e n}{\e J}_{\!\! Y}.
  $$
  It is equally well-known~(cf.~\cite[Example~2.2.5]{Stanley/97})
  that Gaussian polynomials may be expressed in terms of Coxeter
  length functions on symmetric groups:
  $$
  \binom{\e n}{\e J}_{\!\! Y} = \sum_{\substack{w\in\sym_{\e n} \\
      D_\upL(w)\subseteq \e J}} Y^{l(w)}.
  $$
  Equation~\eqref{eq_7} follows in the particular case $I =
  \varnothing$, as $l_{\upL}^{(\e \widetilde{I})^c} = l_{\upL}^S=0$
  and $\lambda_{n,I} = l$.
  
  We now treat the general case $I=\{i_1,\dots,i_r\}_<\subseteq[n-1]$.
  To prove \eqref{eq_7} we argue by induction on~$n$.  The base step
  $n = 0$ is trivial, so suppose that $n > 0$. We may further assume
  that $J = \{ j_1, \dots, j_s \}_< \not = \varnothing$ and we define
  $j := j_1 = \min J$. Our first aim is to derive a recursive formula
  for~$a_{n,I}^J(q)$, using the formula we obtained in the special
  case~$I=\varnothing$. For this purpose we determine what are the
  possible first terms $U_j$ of the flags $\mathbf{U}_J$ we intend to
  count. Then we consider in how many ways each such space $U_j$ can
  be completed to yield a full flag $\mathbf{U}_J$.

  Let $T$ be the set of all $r$-tuples $\mathbf{t} = (t_1,\dots,t_r)
  \in ([j]_0)^r$ such that
  \begin{equation} \label{eq_8}
  \begin{split} 
    & t_1 \leq \dots \leq t_r, \qquad \e \{ t_1, \dots,
    t_r \} \subseteq [\e j]_0 \qquad \text{and} \\
    & \forall \rho\in[r+1]:\; j - (n - i_\rho) - t_{\rho-1} \leq
    t_\rho-t_{\rho-1} \leq i_\rho-i_{\rho-1},
  \end{split}
  \end{equation}
  where $i_0 = t_0 = 0$ and $i_{r+1}:=n, t_{r+1}:=j$. These
  `admissible' tuples encode the possible dimensions of the
  intersections $U_j \cap R_{i_\rho}$ of a $j$-dimensional subspace
  $U_j$ of $V$, non-degenerate with respect to $\bfB$, with the
  members $R_{i_\rho}$ of the flag of radicals associated to
  $\bfB$. Recalling that the underlying field $F$ has
  cardinality $q^{2 \e}$ and applying \eqref{eq_7} for
  $I = \varnothing$, we note that for each $\mathbf{t} \in T$ there
  are precisely
  \begin{equation}\label{eq_9}
  \begin{split}
    A_{n,I}^\mathbf{t}(q) & = \prod_{\rho\in[r+1]}
    a_{(i_\rho-i_{\rho-1}),\varnothing}^{\{t_\rho-t_{\rho-1}\}}
    (q^{-1}) \, q^{2\e(t_\rho-t_{\rho-1})(i_{\rho-1}-t_{\rho-1})} \\
    & = \prod_{\rho \in [r+1]} q^{2\e (t_\rho - t_{\rho -1}) (i_\rho -
      i_{\rho -1}) } \prod_{\rho \in [r+1]} \sum_{\substack{ w \in
        \sym_{\e(i_\rho-i_{\rho-1})} \\ D_\upL(w) \subseteq \e
        \{t_{\rho}-t_{\rho-1}\}}} Y^{l(w)}
  \end{split}
  \end{equation}
  subspaces $U_j$, non-degenerate with respect to $\bfB$, such
  that $\dim( U_j \cap R_{i_\rho} ) = t_\rho$ for all $\rho \in
  [r+1]$. Given such a subspace $U_j$, the number of non-degenerate
  flags $\mathbf{U}_J$ with first term $U_j$ can be described
  inductively, using the notation $J - j = \{ j_2 - j, \dots, j_s - j
  \}$ and $I - \mathbf{t} = \{i_1-t_1,\dots,i_r-t_r\} \cap \N$; it equals
  \begin{equation}\label{eq_10}
    a^{J-j}_{n-j,I-\mathbf{t}}(q) =
    q^{2\e\deg_q\binom{n-j}{J-j}_{\! q}}\sum_{\substack{w\in\sym_{\e(n-j)} \\
      D_\upL(w)\subseteq\e(J-j)}} Y^{\lambda_{n-j,I-\mathbf{t}}(w)}.   
  \end{equation}
  
  For $\mathbf{t} \in T$, apply equations~\eqref{eq_9} and
  \eqref{eq_10} together with the identities
  $$
  \sum_{\rho\in[r+1]}(t_{\rho}-t_{\rho-1})(i_{\rho}-t_{\rho}) =
  j(n-j)-\sum_{\rho\in[r]}t_\rho(i_{\rho+1}-i_\rho-(t_{\rho+1}-t_\rho))
  $$
  and
  $$
  \binom{n}{J}_{\!\! q} = \binom{n}{j}_{\!\! q} \binom{n-j}{J-j}_{\!\! q},\quad
  \deg_q\binom{n}{j}_{\!\! q}=j(n-j)
  $$
  to obtain
  \begin{equation} \label{eq_11}
  \begin{split}
    \alpha_{n,I}^J(q^{-1}) & = a_{n,I}^J(q) / q^{2\e \deg_q
      \binom{n}{J}_{\! q}} \\
    & = q^{-2\e \deg_q \binom{n}{J}_{\! q}} \; \sum_{\mathbf{t} \in T}
    A_{n,I}^\mathbf{t}(q) \;
    a_{n-j,I-\mathbf{t}}^{J-j}(q) \\
    & = \sum_{\mathbf{t} \in T} \; Y^{2\left(\sum_{\rho\in[r]}\e
        t_\rho(\e(i_{\rho+1}-i_\rho)
        - \e(t_{\rho+1} - t_\rho))\right)} \\
    & \quad \cdot \Big( \sum_{\substack{w\in\sym_{\e(n-j)} \\ D_\upL(w)
        \subseteq \e(J-j)}} \!\!\!\! Y^{\lambda_{n-j,I-\mathbf{t}}(w)}
    \Big) \Big( \prod_{\rho\in[r+1]}
    \sum_{\substack{w\in\sym_{\e(i_\rho-i_{\rho-1})} \\
        D_\upL(w) \subseteq \e \{t_{\rho}-t_{\rho-1}\}}} \!\!\!\!
    Y^{l(w)} \Big).
  \end{split}
  \end{equation}
  We are looking to prove that the right hand side of
  equation~\eqref{eq_11} may be written as a sum over the
  elements in the symmetric group~$\sym_{\e n}$ whose left descent set
  is contained in~$\e J$. In the following considerations we shall
  identify permutations $w\in\sym_{\e n}$ (acting on $\{1, \dots, \e
  n\}$ from the right) with the corresponding $\e n \times \e
  n$-permutation matrices (acting on the set of standard row vectors
  by right-multiplication). Observe that for any element $w\in
  \sym_{\e n}$ with $D_\upL(w)\subseteq \e J$ the corresponding
  permutation matrix is ascending on the first~$[\e j]$ rows. Define
  $$
  \mathbf{t} = \mathbf{t}(w) = (t_1,\dots,t_r)
  $$
  by
  $$
  t_\rho := \e^{-1} |\{ \varrho \in[\e j] \mid \varrho^w > \e
  (n-i_\rho)\}| \quad \text{for all $\rho\in[r]$,}
  $$
  and set $t_0 := 0$, $t_{r+1} := j$. Then $\mathbf{t} \in T$,
  as~$n - i_\rho \geq j - t_\rho$ for all~$\rho\in[r+1]$ and thus
  $\mathbf{t}$ satisfies~\eqref{eq_8}. Applying suitable elementary
  column operations to $w$ corresponding to left multiplication by
  elements of the parabolic subgroup $W_{(\e \widetilde{I})^c}$, it is
  easily seen that there are unique elements~$u,v\in \sym_{\e n}$ such
  that
  \begin{enumerate}
  \item[(a)] $w=uv$ and $l(w)=l(u)+l(v)$;
  \item[(b)] for all $\rho \in [r+1]$:
    $$
    \varrho \in [ \e(j-t_\rho)+1, \e(j-t_{\rho-1}) ] \iff
    \varrho^u = \varrho + \e(n-i_\rho);
    $$
  \item[(c)] $v \in W_{(\e \widetilde{I})^c}$, i.e.\ for all $\rho \in
    [r+1]$:
    $$
    \varrho \in [ \e(n-i_\rho)+1, \e(n-i_{\rho-1}) ] \iff
    \varrho^v \in  [\e(n-i_\rho)+1, \e(n-i_{\rho-1}) ],
    $$
    and $D_\upL(v)\subseteq \{\e(n-i_\rho)+\e(t_\rho-t_{\rho-1}) \mid
    \rho\in[r+1]\}$.
  \end{enumerate}
  This is best seen in terms of permutation matrices. We write $\Id_s$
  for the $s\times s$-unit matrix. Then the permutation matrix~$u$
  has the shape

  \begin{equation}\label{eq_12}
  \left(
  \begin{array}{c|c!{\vrule width 1pt}c!{\vrule width 1pt}c|c!{\vrule
        width 1pt}c|c} 
   \Id_{\e(j-t_r)}& & \dots & & & & \\
   \hline
   & & & & & & \\
   \hline
   & & \ddots & & & & \\
   \hline
   & & & \Id_{\e(t_2-t_1)} & & & \\
   \hline
   & & & & & \Id_{\e t_1} & \\
   \hline
   & u_{r+1} & \dots & & u_2 & & u_1
  \end{array}
  \right),
  \end{equation}
  where $u_\rho$ is an $\e (n-j) \times \e ((i_\rho - i_{\rho-1}) -
  (t_\rho - t_{\rho-1}))$-matrix for $\rho \in [r+1]$. The
  permutation matrix~$v$ has the form
  $$
  \left(
  \begin{array}{c!{\vrule width 1pt}c!{\vrule width 1pt}c!{\vrule width 1pt}c}
    v_{r+1} & & & \\
    \hline
    & \ddots & & \\
    \hline
    & & v_2 & \\
    \hline 
    & & & v_1
  \end{array}
  \right),
  $$
  where $v_\rho$ is an $\e(i_\rho - i_{\rho-1}) \times \e(i_\rho -
  i_{\rho-1})$-permutation matrix with at most one descent for $\rho
  \in [r+1]$. We may thus identify~$v$ with
  $$
  (v_1,\dots,v_{r+1})\in \sym_{\e i_1}\times
  \sym_{\e(i_2-i_1)}\times\dots\times\sym_{\e(n-i_r)}
  $$
  and have, by slight abuse of notation, for each $\rho\in[r+1]$,
  \begin{equation}\label{eq_13}
    D_\upL(v_\rho) \subseteq \{\e(t_\rho-t_{\rho-1})\} \cap
    [\e(i_\rho-i_{\rho-1})-1]. 
  \end{equation}

  \begin{remark} 
    The above decomposition $w=uv$ is \emph{not} the one from
    Lemma~\ref{lemma_1}. It is important for our purpose that each
    $v_\rho$ has at most one descent.
  \end{remark}
  
  Note that, by deleting the first~$\e j$ rows and respective columns
  in~\eqref{eq_12}, the element $u$ determines a unique $\e(n-j) \times
  \e(n-j)$-permutation matrix
  $$
  u' := 
  \left(
  \begin{array}{cccc}
  u_{r+1} & \cdots & u_2 & u_1  
  \end{array}
  \right)
  $$
  with descent set~$D_\upL(u') = D_\upL(w) - \e j$. 
  
  As we indicated in Section~\ref{section_coxeter}, it is easy to
  determine the length of a permutation given by a permutation matrix:
  it is simply the number of entries~$0$ in the matrix which are not
  below or to the right of an entry~$1$. Thus
  \begin{align*}
    l(u) & = \sum_{\rho \in [r]} \e t_\rho (\e (i_{\rho+1}-i_\rho) -
    \e (t_{\rho+1}-t_\rho))+l(u'), \\
    l(v) & = \sum_{\rho\in[r+1]} l(v_\rho).
  \end{align*}
  Moreover, the parabolic length of $w$ with respect to $(\e
  \widetilde{I})^c$ is determined by $\mathbf{t}$ and by the parabolic
  length of $u'$ with respect to $(\e(\widetilde{I - \mathbf{t}}))^c$:
  \begin{equation*}
    l_\upL^{(\e \widetilde{I})^c}(w) = \sum_{\rho\in[r]} \e t_\rho (
    \e(i_{\rho+1}-i_\rho) - \e (t_{\rho+1}-t_\rho) ) +
    l_\upL^{(\e(\widetilde{I - \mathbf{t}}))^c}(u').
  \end{equation*}
  This gives  
  \begin{equation} \label{eq_14}
  \begin{split}
    \lambda_{n,I}(w) & = (l(u) + l(v)) + l_\upL^{(\e \widetilde{I})^c}(w) \\
    & = 2 \sum_{\rho\in[r]} \e t_\rho (\e (i_{\rho+1}-i_\rho) - \e
    (t_{\rho+1}-t_\rho) ) \\
    & \quad + \lambda_{\e(n-j),{I-\mathbf{t}}}(u') +
    \sum_{\rho\in[r+1]} l(v_s).
  \end{split}
  \end{equation}
  
  Conversely, any $\mathbf{t} \in T$ and any permutations $v_1, \dots,
  v_{r+1}, u'$ of the appropriate degrees such that \eqref{eq_13}
  holds give rise to a permutation~$w$ satisfying
  $$
  D_\upL(w) \subseteq (D_\upL(u')+\e j)\cup\{\e j\}.
  $$
  Thus \eqref{eq_14} shows that the right hand side of
  \eqref{eq_11} is indeed equal to the right hand side of
  \eqref{eq_7}.
\end{proof}


\section{The orthogonal case}\label{section_orthogonal}

Our aim in this section is to complete the proof of
Theorem~\ref{theorem_A}. We consider non-degenerate quadratic spaces
$\calV = (V,B,f)$, where $V$ is an $n$-dimensional vector space over
the finite field~$F = \F_q$, equipped with a quadratic form $f$, and
$B$ denotes the bilinear form obtained by polarising $f$. So for
all $x,y \in V$,
$$
B(x,y) =
\begin{cases}
  f(x+y) + f(x) + f(y) & \text{if $\cha F = 2$,} \\
  \frac{1}{2} \left( f(x+y) - f(x) - f(y) \right) & \text{if $\cha
    F \ne 2$.}
\end{cases}
$$
If $\cha F \ne 2$, then $B$ is non-degenerate symmetric and, as
$f(x) = B(x,x)$ for all $x \in V$, the quadratic form $f$ can easily
be recovered from $B$. If $\cha F = 2$, then $B$ is alternating,
possibly degenerate and carries less information than $f$.

Non-degenerate quadratic spaces over finite fields have been
classified and can be described up to isomorphism as follows;
cf.~\cite[p.~144]{Artin/57}, \cite[Section~3.3]{Cameron/91}. If $\cha
F \ne 2$, then for any given dimension~$n$ there are two possible
isomorphism types of non-degenerate quadratic spaces $\calV =
(V,B,f)$, namely
\begin{align*}
  \text{for $n$ odd:} && \calV & = \calH_1 \perp \dots \perp
  \calH_m \perp \calA_{1,1}, \\
  && \calV & = \calH_1 \perp \dots \perp \calH_m \perp \calA_{1,-1}, \\
  \text{for $n$ even:} && \calV & = \calH_1 \perp \dots
  \perp \calH_{m-1} \perp \calH_m, \\
  && \calV & = \calH_1 \perp \dots \perp \calH_{m-1} \perp \calA_2,
\end{align*}
where $m=\lfloor \frac{n}{2} \rfloor$, the $\calH_i$ denote hyperbolic
planes, $\calA_{1,1}$ (respectively $\calA_{1,-1}$) stands for an
anisotropic line $\langle x \rangle$ with $f(x) \in (F^*)^2$
(respectively $f(x) \in F^* \setminus (F^*)^2$) and $\calA_2$ is an
anisotropic plane. For the purpose of counting non-degenerate flags in
quadratic spaces~$\calV$ of given odd dimension, there is no
significant difference between the two possible isomorphism types.

We now discuss the case $\cha F = 2$. Then the above list still
provides all isomorphism types of non-degenerate quadratic spaces, but
becomes one term shorter: as every element of $F$ is a square, in
any given odd dimension there is (up to isomorphism) just one
non-degenerate quadratic space. In any given even dimension there are
still two isomorphism types. Note also that non-degenerate quadratic
spaces of odd dimension are defective with $1$-dimensional radical,
whereas non-degenerate quadratic spaces of even dimension are
non-defective.

Returning to the task of proving Theorem~\ref{theorem_A}, we recall
from the introduction that, in the even-dimensional case, we attach a
sign~$\epsilon=1$ or $\epsilon=-1$ to~$\calV$ according to whether the
anisotropic kernel of~$\calV$ is $0$- or $2$-dimensional. More
suggestively, we write~$a^J_{2m+1}(q) := a^J_{\calV}(q)$ if $n = 2m+1$
is odd and, similarly, $a^J_{2m,\epsilon}(q) := a^J_{\calV}(q)$ if~$n
= 2m$ is even and $\calV$ of type~$\epsilon$. We are interested in the
polynomials
\begin{align*}
  \alpha^J_{2m+1}(q^{-1}) & := a^J_{2m+1}(q)/q^{\deg_qa^J_{2m+1}},
  \\
  \alpha^J_{2m,\epsilon}(q^{-1}) & :=
  a^J_{2m,\epsilon}(q)/q^{\deg_qa^J_{2m,\epsilon}}.
\end{align*}

\begin{definition}\label{definition_6}
  Given a family $\bfF=\left(F_J(\bfX)\right)_{J\subseteq[n-1]}$ of
  rational functions with the inversion property~\eqref{eq_IP} we
  define respectively
  \begin{align*}
    \Ig_{2m+1}(q^{-1},\bfX) :=\Ig_{2m+1,\bfF}(q^{-1},\bfX) &:=
    \sum_{J\subseteq [n-1]}\alpha_{2m+1}^J(q^{-1})F_J(\bfX), \\
    \Ig_{2m,\epsilon}(q^{-1},\bfX) :=
    \Ig_{2m,\epsilon,\bfF}(q^{-1},\bfX) & := \sum_{J\subseteq
      [n-1]}\alpha_{2m,\epsilon}^J(q^{-1})F_J(\bfX).
  \end{align*}
\end{definition}

To streamline notation, we will sometimes add in the odd-dimensional
case a superfluous $\epsilon$ to expressions like $a^J_n(q)$,
$\alpha^J_n(q^{-1})$ or $\Ig_n(q^{-1},\bfX)$, thus writing e.g.\ 
$a^J_{n,\epsilon}(q)$, $\alpha^J_{n,\epsilon}(q^{-1})$ or
$\Ig_{n,\epsilon}(q^{-1},\bfX)$, irrespective of the parity of $n$.

We now fix a family of rational functions~$\bfF = (F_J(\bfX))_{J
  \subseteq [n-1]}$ with the inversion property~\eqref{eq_IP}. The
assertion of Theorem~\ref{theorem_A} in the orthogonal case then takes
the following form.

\begin{theorem} \label{theorem_2} 
  The Igusa-type functions satisfy functional equations
  \begin{align*}
    \Ig_{2m+1}(q, \bfX^{-1}) & = (-1)^m q^{m^2+m} \Ig_{2m+1}(q^{-1},\bfX),\\
    \Ig_{2m,\epsilon}(q, \bfX^{-1})& = -\epsilon (-1)^m q^{m^2}
    \Ig_{2m,\epsilon}(q^{-1},\bfX).
  \end{align*}
\end{theorem}

We first give an outline of the proof of Theorem~\ref{theorem_2},
deferring precise definitions for a moment. In
Proposition~\ref{proposition_3} we derive explicit
formulae for the polynomials $\alpha^J_{n,\epsilon}(q^{-1})$ from the
well-known formulae for the orders of the orthogonal groups. A key
observation is that the map $J \mapsto \alpha^J_{n,\epsilon}(q^{-1})$
factors over a `bisecting' map $\phi: \mathcal{P}([n]) \rightarrow
\mathcal{P}([m])$. We are thus led to define, for $G \subseteq [m]$,
$I \in \phi^{-1}(G)$,
\begin{equation*}
 \alpha^{\uparrow G}_{n,\epsilon}(q^{-1}) :=
 \alpha^I_{n,\epsilon}(q^{-1})
\end{equation*}
and
$$
F_{\phi^{-1}(G)}(\bfX):=\sum_{I\in\phi^{-1}(G)}F_I(\bfX)
$$
so that
$$
\Ig_{n,\epsilon}(q^{-1},\bfX) = \sum_{G\subseteq[m]}
\alpha^{\uparrow G}_{n,\epsilon}(q^{-1}) F_{\phi^{-1}(G)}(\bfX).
$$
As we shall see, any subset~$G \subseteq [m]$ induces in a natural
way a composition $C(G) := C(G,m)$ of a non-negative integer $N(G)
\leq m$. For $G,H \subseteq [m]$, we denote by $\|G\|$ the number of
parts of $C(G)$ and by $c_{G,H}$ the number of ways the
composition~$C(H)$ refines a truncation of the composition~$C(G)$. We
then prove the following `inversion equations'.

\begin{proposition}\label{proposition_2}
 \begin{enumerate}
   \item[(i)] For each~$H\subseteq[m]$,
    \begin{equation*}
     F_{\phi^{-1}(H)}(\bfX^{-1}) = (-1)^{n-1 + \|H\|} \sum_{G\subseteq[m]}
     c_{G,H} F_{\phi^{-1}(G)}(\bfX). 
    \end{equation*}
  \item[(ii)] For each~$G\subseteq[m]$,
    \begin{equation*}
    \begin{split}
      \alpha^{\uparrow G}_{2m+1}(q) & = (-1)^m q^{m^2+m} \sum_{H
        \subseteq [m]} (-1)^{\|H\|} c_{G,H} \, \alpha^{\uparrow
        H}_{2m+1}(q^{-1}), \\
      \alpha^{\uparrow G}_{2m,\epsilon}(q) & = \epsilon (-1)^m q^{m^2}
      \sum_{H \subseteq [m]} (-1)^{\|H\|} c_{G,H} \, \alpha^{\uparrow
        H}_{2m,\epsilon} (q^{-1}).
    \end{split}
    \end{equation*} 
 \end{enumerate}
\end{proposition}

Theorem~\ref{theorem_2} is an immediate consequence of
Proposition~\ref{proposition_2}: indeed, in the
odd-dimensional case,
\begin{align*}
  \Ig_{2m+1}(q,\bfX^{-1}) & = \sum_{G\subseteq[m]} \,
  \alpha^{\uparrow G}_{2m+1}(q) F_{\phi^{-1}(G)}(\bfX^{-1}) \\
  & = (-1)^m q^{m^2+m} \sum_{G,H\subseteq[m]} (-1)^{\|H\|}c_{G,H} \,
  \alpha^{\uparrow H}_{2m+1}(q^{-1}) F_{\phi^{-1}(G)}(\bfX^{-1}) \\
  & = (-1)^m q^{m^2+m} \sum_{H\subseteq[m]} \alpha^{\uparrow
    H}_{2m+1}(q^{-1}) F_{\phi^{-1}(H)}(\bfX) \\
  & = (-1)^m q^{m^2+m} \Ig_{2m+1}(q^{-1},\bfX).
\end{align*}
The functional equation for $\Ig_{2m,\epsilon}(q^{-1},\bfX)$ follows
in a similar way. In the remainder of this section we give precise
definitions of the above concepts, and we supply a proof of
Proposition~\ref{proposition_2}.

\begin{definition}[Integer compositions]
  By a \emph{composition} $C$ of a non-negative integer $N$ into
  $\rho$ parts we mean a tuple $(x_1, \dots, x_\rho) \in \N^\rho$ such
  that $N = x_1 + \dots + x_\rho$. 
  
  Given $I = \{ i_1, \dots, i_r \}_< \subseteq [n]$, we define 
  $$
  N(I,n) := \max ( [n]_0 \setminus I ) \quad \text{and} \quad
  \rho := \max \{ \varrho \in [r+1]_0 \mid i_{\varrho-1} < N(I,n) \},
  $$
  where $i_{-1} := -1$, $i_0 := 0$. Then $I$ induces a composition
  $C(I,n)$ of $N(I,n)$ into $\|I\|_n := \rho$ parts, namely
  $$
  C(I,n) := (i_1, i_2 - i_1, \dots, i_{\rho - 1} - i_{\rho - 2},
  N(I,n) - i_{\rho - 1} ).
  $$
  Note that, if $I \subseteq [n-1]$, then $N(I,n) = n$ and $\rho =
  r+1$. The map $I \mapsto C(I,n)$ induces a bijection from
  $\mathcal{P}([n-1])$ onto the set of all compositions of $n$.
  
  We define the \emph{bisecting map}
  $$
  \phi: \mathcal{P}([n]) \rightarrow \mathcal{P}([m])
  $$
  as follows: for $I \subseteq [n]$ with $C(I,n) = (x_1, \dots,
  x_\rho)$ set
  \begin{align*}
    \cut(I) & := \left\lfloor \frac{x_1}{2} \right\rfloor +
    \left\lfloor \frac{x_2}{2}
    \right\rfloor + \dots + \left\lfloor \frac{x_\rho}{2} \right\rfloor, \\
    \phi_0(I) & := \Big\{ \left\lfloor \frac{x_1}{2} \right\rfloor +
    \left\lfloor \frac{x_2}{2} \right\rfloor + \dots + \left\lfloor
      \frac{x_\varrho}{2} \right\rfloor \mid \varrho \in [\rho] \Big\}
    \setminus \{ 0, \cut(I) \} \\
  \end{align*}
  and
  \begin{equation*}
    \phi(I) := \phi_0(I) \cup \left[ \cut(I)+1, m \right].
  \end{equation*}
  Note that $N(\phi(I),m) = \cut(I)$ and $\|\phi(I)\|_m \leq \|I\|_n$.
  Moreover, $\phi$ maps $\mathcal{P}([n-1])$ surjectively onto
  $\mathcal{P}([m])$. For subsets $G \subseteq [m]$ we agree to write
  $N(G) := N(G,m)$ and $\|G\| := \|G\|_m$.
\end{definition}

We now give explicit formulae for the polynomials
$\alpha_{n,\epsilon}^J(q^{-1})$. 

\begin{proposition} \label{proposition_3}
  Let $J \subseteq[n-1]$, $H := \phi(J) \subseteq [m]$, and put $Y :=
  q^{-2}$.
  \begin{enumerate}
  \item[(i)] For $n = 2m+1$ odd,
    \begin{equation*}
    \begin{split}
      \alpha^J_{2m+1}(q^{-1}) & = \binom{N(H)}{\phi_0(J)}_{\!\! Y} \;
      \prod_{i = N(H)+1}^m (1 - Y^i) \\
      & = \binom{m}{H \cup \{N(H)\}}_{\!\! Y} \; (1-Y)^{m - N(H)}.
    \end{split}
    \end{equation*}
  
  \item[(ii)] For $n = 2m$ even,
    \begin{align*}
      \alpha^J_{2m,\epsilon}(q^{-1}) & = \binom{m}{J/2}_{\!\! Y} =
      \binom{m}{H}_{\!\! Y} && \text{if $J
        \subseteq 2\N$,} \\
      \alpha^J_{2m,\epsilon}(q^{-1}) & = \binom{N(H)}{\phi_0(J)}_{\!\!
        Y} \; \frac{\prod_{i =
          N(H)+1}^m (1 - Y^i)}{1 + \epsilon q^{-m}} && \\
      & = \binom{m}{H \cup \{N(H)\}}_{\!\! Y} \; \frac{(1 - Y)^{m -
          N(H)}}{1 + \epsilon q^{-m}} && \text{otherwise.}
    \end{align*}
  \end{enumerate}
\end{proposition}

\begin{proof}
  First we are going to prove the assertions in odd characteristic,
  where the discriminant helps to distinguish isomorphism types of
  quadratic spaces and where we can freely apply Witt's Extension and
  Cancellation Theorem.  Afterwards we explain why the formulae also
  remain true in characteristic $2$.
  
  So first suppose that the underlying field $F = \F_q$ has odd
  characteristic. Recall the formulae for the orders of the respective
  orthogonal groups
  \begin{align*}
    |\textup{O}_{2m+1}(\F_q)| & = 2 q^{m^2} \prod_{i\in[m]}
    (q^{2i}-1) =: p_{2m+1}(q) =: p_{2m+1}, \\
    |\textup{O}_{2m}^\epsilon(\F_q)| & = 2 q^{m^2-m} (q^m-\epsilon)
    \prod_{i\in[m-1]} (q^{2i}-1) =: p_{2m,\epsilon}(q) =:
    p_{2m,\epsilon}
  \end{align*}
  (cf.~\cite[p.~147]{Artin/57}, \cite[Theorem~9.11]{Grove/02}), and
  put
  $$
  p_n^\sharp := p_n^\sharp (q) :=
  \begin{cases}
    p_{2m + 1} & \text{if $n = 2m+1$ odd,} \\
    (q^m + \epsilon) p_{2m,\epsilon} & \text{if $n = 2m$ even}.
  \end{cases}
  $$
  Let $J = \{j_1,\dots,j_s\}_< \subseteq[n-1]$, and put $j_0 := 0$,
  $j_{s+1} := n$. Counting non-degenerate flags $\bfU_J = (U_j)_{j \in
    J}$ of type $J$ in $\calV$ is equivalent to counting (ordered)
  orthogonal decompositions
  \begin{equation}\label{eq_15}
   \calV = \calW_1 \perp \dots \perp \calW_{s+1}
  \end{equation}
  with $\dim \calW_\sigma = k_{\sigma} := j_\sigma - j_{\sigma - 1}$
  for all $\sigma \in [s+1]$. The isomorphism type of such an
  orthogonal decomposition is determined by the discriminants
  $\disc \calW_\sigma \in \F_q^* / (\F_q^*)^2 \cong \{1, -1\}$ of the
  non-degenerate spaces $\calW_\sigma$, $\sigma \in [s+1]$.
  
  Let $\eta \in \{1, -1\}$, according to whether $-1$ is a square in
  $\F_q$ or not. At this point it is advantageous to assign, also to
  an odd-dimensional non-degenerate quadratic space $\calW$ a sign
  $\epsilon(\calW) \in \{ 1, -1 \}$, namely the discriminant of the
  (one-dimensional) anisotropic kernel of $\calW$. We then have
  $\disc \calW  = \epsilon(\calW) \eta^{\lfloor \dim \calW  / 2 \rfloor}$
  for any non-degenerate quadratic space $\calW$, irrespective of the
  parity of $\dim \calW$.
  
  Thus the isomorphism type of an orthogonal decomposition of the
  form~\eqref{eq_15} can be encoded in a tuple
  $\boldsymbol{\epsilon} = (\epsilon_1, \dots, \epsilon_{s+1}) \in \{
  1,-1 \}^{s+1}$ such that $\calW_\sigma$ is of type $\epsilon_\sigma$
  for all $\sigma \in [s+1]$. Moreover, the tuples
  $\boldsymbol{\epsilon}$ which arise in this way are precisely the
  elements of $E := E(\calV) := \{ \boldsymbol{\epsilon} \mid
  \epsilon_1 \cdots \epsilon_{s+1} = \eta^{m - N(\phi(J))} \epsilon
  \}$, and Witt's Extension and Cancellation Theorem implies that the
  number of ordered orthogonal decompositions of isomorphism type
  $\boldsymbol{\epsilon} \in E$ equals
  $$
  \frac{|\textup{O}_{n}^\epsilon(\F_q)|}{|\prod_{\sigma = 1}^{s+1}
    \textup{O}_{k_\sigma}^{\epsilon_\sigma}(\F_q)|} =
  \frac{p_{n,\epsilon}}{\prod_{\sigma = 1}^{s+1} p_{k_\sigma,
      \epsilon_\sigma}};
  $$
  cf.~\cite[p.~147f]{Artin/57}. Setting
  $$
  \mathcal{E}(J) := \{ \sigma \in [s+1] \mid k_\sigma \equiv 0 \mod
  2 \}
  $$
  we thus obtain
  \begin{equation*}
    a_{n,\epsilon}^J(q) = p_{n,\epsilon} \sum_{\boldsymbol{\epsilon}
      \in E} \left( \prod_{\sigma = 1}^{s+1} p_{k_\sigma,
        \epsilon_\sigma} \right)^{-1}
     = \frac{p_{n,\epsilon}}{ \prod_{\sigma = 1}^{s+1}
      p_{k_\sigma}^\sharp } \; \sum_{\boldsymbol{\epsilon} \in E}
    \prod_{\sigma \in \mathcal{E}(J)} (q^{k_\sigma / 2} +
    \epsilon_\sigma).
  \end{equation*}
  Note that
  \begin{align*}
    \sum_{\boldsymbol{\epsilon} \in E} \prod_{\sigma \in
      \mathcal{E}(J)} (q^{k_\sigma / 2} + \epsilon_\sigma) & = 2^s
    \prod_{\sigma \in \mathcal{E}(J)} q^{k_\sigma / 2} && \text{if
      $\mathcal{E}(J) \not = [s+1]$,} \\
    \sum_{\boldsymbol{\epsilon} \in E} \prod_{\sigma \in
      \mathcal{E}(J)} (q^{k_\sigma / 2} + \epsilon_\sigma) & = 2^s
    (q^{n / 2} + \epsilon) && \text{if $\mathcal{E}(J) = [s+1]$.}
 \end{align*}  
 From this the claim follows for $\cha F \ne 2$.
 
 Before turning our attention to the case $\cha F = 2$, we record a
 set of formulae for later use. Let $j \in [n-1]$ and $\delta \in
 \{1,-1\}$. If $j = 2h+1$ is odd, let $a_\calV^j(q)$ denote the number
 of non-degenerate $j$-dimensional subspaces in $\calV$. If $j = 2h$
 is even, let $a_\calV^{j,\delta}(q)$ denote the number of
 non-degenerate $j$-dimensional subspaces of type $\delta$ in $\calV$.
 According to whether $\calV$ is odd- or even-dimensional, we also
 write $a_{2m+1}^{2h+1}(q)$, $a_{2m,\epsilon}^{2h+1}(q)$ in the former
 and $a_{2m+1}^{2h,\delta}(q)$, $a_{2m,\epsilon}^{2h,\delta}(q)$ in
 the latter case. Our calculations above, based on Witt's Extension
 and Cancellation Theorem, show in particular that, if $\cha F \ne 2$,
 \begin{align}
   a_{2m+1}^{2h+1}(q) & =
   \frac{|\textup{O}_{2m+1}(\F_q)|}{|\textup{O}_{2h+1}(\F_q)|
     |\textup{O}_{2m-2h}^+(\F_q)|} +
   \frac{|\textup{O}_{2m+1}(\F_q)|}{|\textup{O}_{2h+1}(\F_q)|
     |\textup{O}_{2m-2h}^-(\F_q)|} \label{eq_16} \\
   & = \frac{2 \, q^{m-h} \, p_{2m+1}}{p_{2h+1} \, p^\sharp_{2m-2h}},
   \notag \\
   a_{2m,\epsilon}^{2h+1}(q) & = 2 \cdot
   \frac{|\textup{O}_{2m}^\epsilon (\F_q)|}{|\textup{O}_{2h+1}(\F_q)|
     |\textup{O}_{2m-2h-1}(\F_q)|}
   = \frac{2 \, p_{2m,\epsilon}}{p_{2h+1} \, p_{2m-2h-1}}, \label{eq_17} \\
   a_{2m+1}^{2h,\delta}(q) & = \frac{|\textup{O}_{2m+1}
     (\F_q)|}{|\textup{O}_{2h}^\delta (\F_q)|
     |\textup{O}_{2m-2h+1}(\F_q)|} = \frac{p_{2m+1}}{p_{2h,\delta} \,
     p_{2m-2h+1}}, \label{eq_18} \\
   a_{2m,\epsilon}^{2h,\delta}(q) & = \frac{|\textup{O}_{2m}^\epsilon
     (\F_q)|}{|\textup{O}_{2h}^\delta (\F_q)|
     |\textup{O}_{2m-2h}^{\delta \epsilon} (\F_q)|} =
   \frac{p_{2m,\epsilon}}{p_{2h,\delta} \, p_{2m-2h,\delta \epsilon}}.
   \label{eq_19}
 \end{align}
 Below we will show that, in fact, also in characteristic $2$ one
 obtains the same polynomials $a_\calV^j(q)$ and
 $a_\calV^{j,\delta}(q)$.  Thus, by induction, the formulae for
 $a_\calV^J(q)$ and $\alpha_\calV^J(q^{-1})$, which we initially
 derived only under the extra assumption $\cha F \ne 2$, also remain
 valid in characteristic $2$.
 
 So suppose that $\cha F = 2$, and let $j \in [n-1]$, $\delta \in
 \{1,-1\}$. Write $j = 2h+1$, if $j$ is odd, and $j = 2h$, if $j$ is
 even. The orders of the respective orthogonal groups are now
 \begin{equation*}
  |\textup{O}_{2m+1}(\F_q)| =  \frac{p_{2m+1}}{2}, \qquad
   |\textup{O}_{2m}^\epsilon(\F_q)| = p_{2m,\epsilon}, 
 \end{equation*}
 and Witt's Extension and Cancellation Theorem still applies to
 non-defective subspaces; cf.~\cite[Theorems~3.12 and~14.48]{Grove/02}
 and \cite[Theorem~3.15]{Cameron/91}. Therefore we immediately obtain
 the counterparts of \eqref{eq_18} and \eqref{eq_19},
 \begin{align*}
   a_{2m+1}^{2h,\delta}(q) & = \frac{|\textup{O}_{2m+1}
     (\F_q)|}{|\textup{O}_{2h}^\delta (\F_q)|
     |\textup{O}_{2m-2h+1}(\F_q)|} = \frac{p_{2m+1}}{p_{2h,\delta} \,
     p_{2m-2h+1}}, \\
   a_{2m,\epsilon}^{2h,\delta}(q) & = \frac{|\textup{O}_{2m}^\epsilon
     (\F_q)|}{|\textup{O}_{2h}^\delta (\F_q)|
     |\textup{O}_{2m-2h}^{\delta \epsilon} (\F_q)|} =
   \frac{p_{2m,\epsilon}}{p_{2h,\delta} \, p_{2m-2h,\delta \epsilon}}.
 \end{align*}
 Next we suppose that $n = 2m+1$ is odd and compute
 $a_{2m+1}^{2h+1}(q)$. If $h = 0$, then we are to count anisotropic
 lines in $\calV$. It is well-known that the polar space associated to
 $\calV$ has $(q^{2m} - 1)/(q-1)$ points, each corresponding to an
 isotropic line; cf.~\cite[Theorem~3.13]{Cameron/91}. So we deduce
 that
 $$
 a_{2m+1}^{1}(q) = \frac{q^{2m+1} - 1}{q-1} - \frac{q^{2m} -
   1}{q-1} = \frac{2 \, q^m \, p_{2m+1}}{p_{1} \, p^\sharp_{2m}}.
 $$
 In general, choosing a $(2h+1)$-dimensional non-degenerate
 subspace $U$ in $\calV$ can be split into two parts: first pick a
 $2h$-dimensional non-degenerate (hence non-defective) subspace $U_0$
 of type $1$, then complement your choice by picking an anisotropic
 line $A$ in $U_0^\perp$ to obtain $U = U_0 + A$. Applying Witt's
 Extension and Cancellation Theorem, we obtain the counterpart of
 \eqref{eq_16},
 \begin{equation*}
 \begin{split}
   a_{2m+1}^{2h+1}(q) & = \frac{a_{2m+1}^{2h,1}(q) \;
     a_{2m-2h+1}^{1}(q)}{a_{2h+1}^{2h,1}(q)} \\
   & = \frac{ ( p_{2m+1} \cdot 2 \, q^{m-h} \, p_{2m-2h+1} ) /
     (p_{2h,1} \, p_{2m-2h+1} \cdot p_{1} \, p^\sharp_{2m-2h} ) }{
     (p_{2h+1}) / (p_{2h,1} \, p_{1})} \\
   & = \frac{2 \, q^{m-h} \, p_{2m+1}}{p_{2h+1} \, p^\sharp_{2m-2h}}.
 \end{split}
 \end{equation*}
 A similar computation yields the counterpart of \eqref{eq_17}.
\end{proof}

\begin{definition}[Refinements of compositions]\label{definition_8}
  Let $C_1 = (x_1, \dots, x_{\kappa})$ and $C_2 = (y_1, \dots,
  y_{\lambda})$ be compositions. A \emph{refinement of a truncation of
    $C_1$ by $C_2$} is a triple $(C_1,C_2,\boldsymbol{\xi})$ such that
  $\boldsymbol{\xi}=(\xi_1,\dots,\xi_{\kappa}) \in
  {[\lambda]_0}^{\kappa}$ satisfies
  $$
  \xi_1 \leq \dots \leq \xi_{\kappa} = \lambda \quad \text{and}
  \quad \forall i \in [\kappa]: \; y_{\xi_{i-1} + 1}+ \dots +y_{\xi_i}
  \leq x_i,
  $$
  where $\xi_0 := 0$. By slight abuse of terminology, we also call
  the $\kappa$-tuple $\boldsymbol{\xi}$ a refinement of a truncation
  of $C_1$ by $C_2$. For $G,H \subseteq [m]$, the number of
  refinements of truncations of $C(G)$ by $C(H)$ is denoted
  by~$c_{G,H} := c_{G,H}^{(m)}$.
  
  Let $I, J \subseteq [n-1]$ such that $I \subseteq J$, and put $G :=
  \phi(I)$, $H := \phi(J)$. Clearly, $C(J,n)$ can be regarded as a
  refinement of $C(I,n)$. Applying the bisecting map, we obtain a
  refinement of a truncation of $C(G)$ by $C(H)$ as follows.
  
  The sets $[n-1] \setminus I$ and $[n-1] \setminus J$ decompose
  uniquely into disjoint unions
  $$
  [n-1] \setminus I = \mathcal{I}_{I,1} \disjunion \dots \disjunion
  \mathcal{I}_{I,\|G\|}, \qquad [n-1] \setminus J = \mathcal{I}_{J,1}
  \disjunion \dots \disjunion \mathcal{I}_{J,\|H\|}
  $$
  of intervals $\mathcal{I}_{I,i}$ (respectively $\mathcal{I}_{J,j}$)
  of natural numbers with $\max \mathcal{I}_{I,i} < \min
  \mathcal{I}_{I,i+1}$ (respectively $\max \mathcal{I}_{J,j} < \min
  \mathcal{I}_{J,j+1}$) for all admissible values of $i$ (respectively
  $j$).
  
  The \emph{refinement of a truncation of $G$ by $H$ induced from~$I
    \subseteq J$} is the $\|G\|$-tuple $\boldsymbol{\xi}(I,J) =
  (\xi_1,\dots,\xi_{\|G\|})$ defined by $\xi_{\|G\|}:=\|H\|$ and
  \begin{equation*}
    \forall i \in [\|G\|] : \; \mathcal{I}_{J,\xi_{i-1} + 1}
    \disjunion \dots \disjunion \mathcal{I}_{J,\xi_i}  \subseteq
    \mathcal{I}_{I,i}, 
  \end{equation*}
  where $\xi_0 := 0$. We remark that, starting from $G,H \subseteq
  [m]$, every refinement of a truncation of $C(G)$ by $C(H)$ is
  induced by suitable $I,J \subseteq [n-1]$ with $I \subseteq J$.
\end{definition}

We illustrate these notions by an example. 

\begin{example}
  Set $n=11$ so that $m=5$. The subsets $G=\{1,3,4\}$, $H=\{2,4,5\}
  \subseteq [m]$ induce compositions $C(G)=(1,2,1,1)$ of $N(G)=5$ and
  $C(H) = (2,1)$ of $N(H)=3$, respectively. Note that~$\|G\|=4$
  and~$\|H\|=2$.  Among the seven `truncations' of $(1,2,1,1)$ to
  `pre-compositions' of $3$,
  
  \medskip 
  
  \centerline{
  \begin{tabular}{ccccc}
      $(1,2,0,0)$, &$(1,1,1,0)$, &$(1,1,0,1)$, &$(1,0,1,1)$,
      &$(0,1,1,1)$, \\ 
      $(0,2,0,1)$, &$(0,2,1,0)$, &             &             &
  \end{tabular}
  }
  
  \medskip

  \noindent only the $2 = c_{G,H}$ last ones yield the composition
  $C(H) = (2,1)$. They are encoded in the tuples
  $\boldsymbol{\xi}=(0,1,1,2)$ and $\boldsymbol{\xi}=(0,1,2,2)$,
  respectively. 
  
  Define subsets
  \begin{equation*}
    I := \{2,7,9\}, \qquad
    J_1 := \{1,2,7,8,9\}, \qquad
    J_2 := \{1,2,3,7,9,10\}
  \end{equation*}
  of $[n-1] = [10]$. The set $I$ induces the composition $C(I,n) =
  (2,5,2,2)$ of $N(I,n) = 11$. Thus $\cut(I) = 1+2+1+1 = 5$ and
  $\phi(I) = G$. Similarly, $\cut(J_1) = \cut(J_2) = 2+1 =3$ and
  $\phi(J_1)=\phi(J_2) = H$. We have
  $\boldsymbol{\xi}(I,J_1)=(0,1,1,2)$ and
  $\boldsymbol{\xi}(I,J_2)=(0,1,2,2)$.
\end{example}

We are now ready to prove Proposition~\ref{proposition_2}.

\begin{proof}[Proof of Proposition~\textup{\ref{proposition_2}~(i)}] 
  
  Let $H \subseteq [m]$. From the definition of
  $F_{\phi^{-1}(H)}(\bfX)$ and the fact that $\bfF$ has the inversion
  property~\eqref{eq_IP} we obtain
  $$
  F_{\phi^{-1}(H)}(\bfX^{-1}) = \sum_{J\in\phi^{-1}(H)} (-1)^{|J|}
  \sum_{I \subseteq J} F_I(\bfX).
  $$
  Thus it is enough to show that for $I \subseteq [n-1]$ with
  $\phi(I)=G$,
  $$
  (-1)^{n-1+\|H\|} \sum_{\substack{J\in\phi^{-1}(H) \\ I \subseteq
      J}} (-1)^{|J|} = c_{G,H}.
  $$
  This is certainly the case if $c_{G,H} = 0$, as then the sum on
  the left hand side is empty.  Now suppose that $c_{G,H} \not =
  0$ and fix a refinement~$\boldsymbol{\xi}$ of a truncation
  of~$G$ by~$H$; put $\xi_0 := 0$. It suffices to show that
  \begin{equation}\label{eq_20}
    (-1)^{n-1+\|H\|} \sum_{\substack{J\in\phi^{-1}(H) \\ I \subseteq
    J,\;\; \boldsymbol{\xi}(I,J) = \boldsymbol{\xi}}} (-1)^{|J|} = 1.
  \end{equation} 
  
  Decompose $[n-1] \setminus I = \mathcal{I}_{I,1} \disjunion \dots
  \disjunion \mathcal{I}_{I,\|G\|}$ into a disjoint union of intervals
  $\mathcal{I}_{I,i}$ as in Definition~\ref{definition_8},
  and write $C(H) = (y_1, \dots, y_{\|H\|})$. We claim that
  \begin{multline}\label{eq_21}
    \sum_{\substack{J\in\phi^{-1}(H) \\ I \subseteq J,\;\;
        \boldsymbol{\xi}(I,J) = \boldsymbol{\xi}}} (-1)^{|J| - |I|} =
    \prod_{i=1}^{\|G\|} \; \sum_{k = 0}^{\xi_i-\xi_{i-1}}
    \binom{|\mathcal{I}_{I,i}| - \sum_{j = \xi_{i-1} + 1}^{\xi_i}
      (2y_j - 1) - k + 1}{\xi_i - \xi_{i-1}} \\ \cdot \binom{\xi_i -
      \xi_{i-1}}{k} (-1)^{|\mathcal{I}_{I,i}| - \sum_{j = \xi_{i-1} +
        1}^{\xi_i} (2y_j - 1) - k}.
  \end{multline}
   Indeed, specifying $J \in\phi^{-1}(H)$ with $I \subseteq J$ and
  $\boldsymbol{\xi}(I,J) = \boldsymbol{\xi}$ is equivalent to the
  following task: for each $i \in [\|G\|]$ choose $k_i \in
  [\xi_i-\xi_{i-1}]_0$ and single out a disjoint union
  $\mathcal{I}_{J,\xi_{i-1} + 1} \disjunion \dots \disjunion
  \mathcal{I}_{J,\xi_i} \subseteq \mathcal{I}_{I,i}$ of intervals
  $\mathcal{I}_{J,j}$ such that
  \begin{enumerate}
  \item[(a)] $\max \mathcal{I}_{J,j} < \min \mathcal{I}_{J,j+1}$ for
    all admissible values of $j$, 
  \item[(b)] $|\mathcal{I}_{J,j}| = 2 y_j$ for exactly $k_i$ values
    of $j$ and $|\mathcal{I}_{J,j}| = 2 y_j - 1$ for the remaining
    values of $j$.
  \end{enumerate}
  Moreover, the cardinality of the set $J$ corresponding to such a
  choice of $k_i$ and such a choice of intervals $\mathcal{I}_{J,j}
  \subseteq \mathcal{I}_{I,i}$ is
  $$
  |I| + \sum_{i = 1}^{\|G\|} \left( |\mathcal{I}_{I,i}| - \sum_{j =
      \xi_{i-1} + 1}^{\xi_i} (2 y_j -1) - k_i \right).
  $$
  
  As $n-1 + \|H\| = |I| + \sum_{i = 1}^{\|G\|} |\mathcal{I}_{I,i}| +
  \sum_{i = 1}^{\|G\|} (\xi_i - \xi_{i-1})$,
  equation~\eqref{eq_21} implies that the left hand side
  of~\eqref{eq_20} is equal to
  $$
  \prod_{i=1}^{\|G\|} \sum_{k=0}^{\xi_i-\xi_{i-1}}
  \binom{|\mathcal{I}_{I,i}| - \sum_{j = \xi_{i-1} + 1}^{\xi_i} (2y_j
    -1) + 1 - k}{\xi_i - \xi_{i-1}} \binom{\xi_i -
    \xi_{i-1}}{k}(-1)^k.
  $$
  This does indeed equal $1$, because for any positive integers $M
  \leq N$,
  $$
  \sum_{k=0}^M \binom{N-k}{M} \binom{M}{k} (-1)^k = 1
  $$
  (cf.~\cite[p.~169, (5.25)]{GrahamKnuthPatashnik/98}) and hence
  each of the $\|G\|$ factors already equals $1$.
\end{proof}

\begin{proof}[Proof of Proposition~\textup{\ref{proposition_2}~(ii)}
  for $n = 2m+1$ odd]
  
  For $G \subseteq [m]$, we are looking to prove
  \begin{equation}\label{eq_22}
    \alpha_{2m+1}^{\uparrow G}(q) = (-1)^m q^{m^2+m} \sum_{H
    \subseteq [m]} (-1)^{\|H\|} c_{G,H} \; \alpha_{2m+1}^{\uparrow H}
    (q^{-1}).
  \end{equation}
  First we deal with the case $N(G) < m$, i.e.\ $m \in G$. Writing $G'
  := G \setminus \{m\}$ and $Y := q^{-2}$, we see from
  Proposition~\ref{proposition_3}~(i) that in this
  case
  $$
  \alpha_{2m+1}^{\uparrow G}(q^{-1}) = (1 - Y^m)
  \alpha_{2m-1}^{\uparrow G'} (q^{-1}).
  $$
  If $H \subseteq [m]$ with $c_{G,H} \not = 0$, then $N(H) \leq
  N(G) < m$, hence $m \in H$, hence we obtain $\|H\| = \|H'\|_{m-1}$
  and $c_{G,H} = c_{G',H'}^{(m-1)}$ for $H' := H \setminus \{m\}$.
  With these observations \eqref{eq_22} follows
  by induction:
  \begin{align*}
    \alpha_{2m+1}^{\uparrow G}(q) & = (1 - Y^{-m}) (-1)^{m-1}
    q^{(m-1)^2+(m-1)} \\
    & \qquad \cdot \sum_{H' \subseteq [m-1]} (-1)^{\|H'\|_{m-1}}
    c_{G',H'}^{(m-1)} \; \alpha_{2m-1}^{\uparrow H'} (q^{-1}) \\
    & = (-1)^m q^{m^2+m} \sum_{H' \subseteq [m-1]} (-1)^{\|H'\|_{m-1}}
    c_{G',H'}^{(m-1)} \; (1 - Y^m) \alpha_{2m-1}^{\uparrow H'} (q^{-1}) \\
    & = (-1)^m q^{m^2+m} \sum_{H \subseteq [m]} (-1)^{\|H\|} c_{G,H}
    \; \alpha_{2m-1}^{\uparrow H} (q^{-1}).
  \end{align*}
  
  It remains to consider the case $N(G) = m$, i.e.\ $G \subseteq
  [m-1]$. Again set $Y := q^{-2}$, and write $C(G) = (x_1, \dots,
  x_{k+1})$. Proposition~\ref{proposition_3}~(i)
  shows that in this case
  $$
  \alpha_{2m+1}^{\uparrow G}(q^{-1}) = \binom{m}{G}_{\!\! Y},
  $$
  in particular, as $\deg_Y \binom{m}{G}_{\! Y} = \binom{m+1}{2} -
  \sum_{\kappa \in [k+1]} \binom{x_\kappa + 1}{2}$,
  $$
  \alpha_{2m+1}^{\uparrow G}(q) =
  \alpha_{2m+1}^{\uparrow G}(q^{-1}) \; Y^{- \binom{m+1}{2} +
    \sum_{\kappa \in [k+1]} \binom{x_\kappa + 1}{2}} .
  $$
  We shall show below that 
  \begin{equation} \label{eq_23}
    \alpha_{2m+1}^{\uparrow G}(q^{-1}) \; Y^{\sum_{\kappa \in [k+1]}
    \binom{x_\kappa + 1}{2}} = \sum_{H \subseteq [m]}
    (-1)^{m + \|H\|} c_{G,H} \; \alpha_{2m+1}^{\uparrow H}(q^{-1}).
  \end{equation}
  From these equations \eqref{eq_22} follows readily.
  
  It remains to prove \eqref{eq_23}.  For this we need
  the following formulae.
  \begin{enumerate}
  \item[(i)] For all $i \in \N_0$: $Y^{\binom{i+1}{2}} = \sum_{j=0}^i
    \binom{i}{[i-j, i-1]}_{\! Y} (Y-1)^j Y^{\binom{i-j}{2}}$.
  \item[(ii)] For all $i \in \N_0$: $Y^{\binom{i+1}{2}} =
    \sum_{I\subseteq[i]}\binom{i+1}{I}_{\! Y}(-1)^{i-|I|}$.
  \end{enumerate}
  Part (i) is easily proved inductively (see the end of this proof),
  part (ii) is a well-known fact about Gaussian polynomials. With the
  formulae (i), (ii) at our disposal, the left hand side of
  \eqref{eq_23} can be written as
  \begin{equation}\label{eq_24}
  \begin{split}
    \alpha_{2m+1}^{\uparrow G} & (q^{-1}) \; Y^{\sum_{\kappa \in
        [k+1]} \binom{x_\kappa + 1}{2}} \\
    & = \binom{m}{G}_{\!\! Y} \prod_{\kappa \in [k+1]} Y^{\binom{x_\kappa +
        1}{2}} \\
    & = \binom{m}{G}_{\!\! Y} \prod_{\kappa \in [k+1]} \left(
      \sum_{j=0}^{x_\kappa} \binom{x_\kappa}{[ x_\kappa-j, x_\kappa-1
        ]}_{\!\! Y} \left(Y - 1 \right)^j
      Y^{\binom{x_\kappa - j}{2}} \right) \\
    & = \binom{m}{G}_{\!\! Y} \prod_{\kappa \in [k+1]} \left(
      \sum_{j=0}^{x_\kappa} \binom{x_\kappa}{[x_\kappa-j,
        x_\kappa-1 ]}_{\!\! Y} \left( 1 - Y \right)^j \right. \\
    & \qquad \qquad \qquad \qquad \qquad \cdot \left. \sum_{K
        \subseteq [x_\kappa-j-1]} \binom{x_\kappa - j}{K}_{\!\! Y}
      (-1)^{x_\kappa + |K| + \delta(\text{`} j \not = x_\kappa
        \text{'})} \right),
  \end{split}
  \end{equation}
  where the Kronecker-delta $\delta(\text{`} j \not = x_\kappa
  \text{'}) \in \{1,0\}$ reflects whether or not the inequality $j
  \not = x_\kappa$ holds. On the other hand, setting
  $$
  \Xi := \bigcup_{H \subseteq [m]} \{ (H,\boldsymbol{\xi}) \mid
  \boldsymbol{\xi} \text{ a refinement of a truncation of $G$ by $H$}
  \},
  $$
  the right hand side of~\eqref{eq_23} can be
  written as
  \begin{equation}\label{eq_25}
  \sum_{H \subseteq [m]} (-1)^{m + \|H\|} c_{G,H} \;
  \alpha_{2m+1}^{\uparrow H}(q^{-1}) = \sum_{(H,\boldsymbol{\xi}) \in
    \Xi} (-1)^{m + \|H\|} \alpha_{2m+1}^{\uparrow H}(q^{-1}).
  \end{equation}
  Now we explain why the last sum is indeed equal to the right hand
  side of \eqref{eq_24}. Choosing an element
  $(H,\boldsymbol{\xi}) \in \Xi$, so that $\boldsymbol{\xi} = (\xi_1,
  \dots, \xi_{k+1})$ is a refinement of a truncation of $C(G) = (x_1,
  \dots, x_{k+1})$ by $C(H) = (y_1, \dots, y_\lambda)$, is the same as
  fixing for each $\kappa \in [k+1]$ a truncation length $j_\kappa \in
  [x_\kappa]_0$ and a subset $K_\kappa \subseteq [x_\kappa - j_\kappa
  - 1]$, corresponding to a composition $(y_{\xi_{\kappa - 1} + 1},
  \dots, y_{\xi_\kappa})$ of $x_\kappa - j_\kappa$. Moreover, the
  summands attached to the data $(H,\boldsymbol{\xi})$ in
  \eqref{eq_25} and $(j_\kappa,K_\kappa)_{\kappa \in [k+1]}$
  in \eqref{eq_24} respectively agree:
  \begin{equation*}
    (-1)^{m + \|H\|} = (-1)^{\sum_{\kappa \in [k+1]} x_\kappa +
      \sum_{\kappa \in [k+1]} (|K_\kappa| + \delta(\text{`} j \not =
      x_\kappa \text{'}) )} 
  \end{equation*}
  and
  \begin{align*}  
    \alpha_{2m+1}^{\uparrow H}(q^{-1}) & = \frac{(1 - Y^m) (1 -
      Y^{m-1}) \cdots (1-Y)}{\prod_{\iota \in [\lambda]} (1 -
      Y^{y_\iota}) (1 - Y^{y_\iota - 1}) \cdots (1 - Y)} \\
    & = \frac{(1 - Y^m) (1 - Y^{m-1}) \cdots (1-Y)}{\prod_{\kappa \in
        [k+1]} \prod_{\iota = \xi_{\kappa - 1} + 1}^{\xi_\kappa} (1 -
      Y^{y_\iota}) (1 - Y^{y_\iota - 1}) \cdots (1 - Y)} \\
    & = \frac{(1 - Y^m) (1 - Y^{m-1}) \cdots (1-Y)}{\prod_{\kappa \in
        [k+1]} (1 - Y^{x_\kappa - j_\kappa}) (1 - Y^{x_\kappa -
        j_\kappa - 1}) \cdots (1 - Y)} \prod_{\kappa \in [k+1]} \!
    \binom{x_\kappa - j_\kappa}{K_\kappa}_{\!\! Y} \\
    & = \frac{(1 - Y^m) (1 - Y^{m-1}) \cdots (1-Y)}{\prod_{\kappa \in
        [k+1]} (1 - Y^{x_\kappa}) (1 - Y^{x_\kappa - 1}) \cdots (1
      - Y)}   \\
    & \qquad \cdot \prod_{\kappa \in [k+1]}
    \binom{x_\kappa}{[x_{\kappa - j_\kappa}, x_\kappa - 1]}_{\!\! Y} (1 -
    Y)^{j_\kappa} \binom{x_\kappa -
      j_\kappa}{K_\kappa}_{\!\! Y} \\
    & = \binom{m}{G}_{\!\! Y} \prod_{\kappa \in [k+1]}
    \binom{x_\kappa}{[x_{\kappa - j_\kappa}, x_\kappa - 1]}_{\!\! Y} (1 -
    Y)^{j_\kappa} \binom{x_\kappa - j_\kappa}{K_\kappa}_{\!\! Y}.
  \end{align*}
  This finishes the proof of \eqref{eq_23}. For later
  use we record
  \begin{equation} \label{eq_26}
    \alpha_{2m+1}^{\uparrow G}(q^{-1}) \; Y^{\sum_{\kappa \in [k+1]}
    \binom{x_\kappa}{2}} = \sum_{\substack{H \subseteq [m] \\ N(H) = m}}
    (-1)^{m + \|H\|} c_{G,H} \; \alpha_{2m+1}^{\uparrow H}(q^{-1}).
  \end{equation}
  Indeed, summing only over those $H \subseteq [m]$ such that $N(H) =
  m$ is achieved by setting persistently $j = j_\kappa = 0$ in the
  above formulae. Clearly, under the restriction $j = 0$ the term in
  the third line of~\eqref{eq_24} reduces to the left hand
  side of \eqref{eq_26}.
  
  Finally, we supply the proof of the formulae (i) above. We argue by
  induction on $i \in \N_0$. For $i=0$ we have
  $$
  Y^{\binom{1}{2}} = 1 = \binom{0}{\varnothing}_{\!\! Y} (Y-1)^0
  Y^{\binom{0}{2}},
  $$
  and for $i>0$ we obtain, by induction,
  \begin{align*}
    \sum_{j=0}^i & \binom{i}{[i-j,i-1]}_{\!\! Y} (Y-1)^j Y^{\binom{i-j}{2}}
    \\
    & = Y^{\binom{i}{2}} + \sum_{j \in [i]} \binom{i}{[i-j,i-1]}_{\!\! Y} (Y-1)^j
    Y^{\binom{i-j}{2}}\\
    & = Y^{\binom{i}{2}} + \sum_{j \in [i]} \binom{i-1}{[i-j,i-2]}_{\!\! Y}
    \binom{i}{i-1}_{\!\! Y} (Y-1)^j
    Y^{\binom{i-j}{2}} \\
    & = Y^{\binom{i}{2}} + (Y-1) \binom{i}{i-1}_{\!\! Y} \sum_{j \in [i]}
    \binom{i-1}{[i-j,i-2]}_{\!\! Y} (Y-1)^{j-1}
    Y^{\binom{i-j}{2}} \\
    & = Y^{\binom{i}{2}} + (Y-1) \frac{(Y^i -1)}{(Y-1)}
    Y^{\binom{i}{2}} = Y^{i+\binom{i}{2}} = Y^{\binom{i+1}{2}}.
  \end{align*}
\end{proof}

\begin{proof}[Proof of Proposition~\textup{\ref{proposition_2}~(ii)}
    for $n = 2m$ even] 
  
  For $G \subseteq [m]$ we are looking to prove
  \begin{equation} \label{eq_27}
  \alpha_{2m,\epsilon}^{\uparrow G}(q) = \epsilon (-1)^m q^{m^2}
  \sum_{H \subseteq [m]} (-1)^{\|H\|} c_{G,H} \;
  \alpha_{2m,\epsilon}^{\uparrow H}(q^{-1}).
  \end{equation}
  Again by an inductive argument, analogous to the case $n = 2m+1$, we
  may assume that in fact $N(G)=m$, i.e.\ $G \subseteq [m-1]$.  Write
  $Y := q^{-2}$ and $C(G) = (x_1, \dots, x_{k+1})$.
  Proposition~\ref{proposition_3}~(ii) shows that
  in this case
  $$
  \alpha_{2m,\epsilon}^{\uparrow G}(q^{-1}) = \alpha_{2m +
    1}^{\uparrow G}(q^{-1}) = \binom{m}{G}_{\!\! Y},
  $$
  in particular, as $\deg_Y \binom{m}{G}_{\! Y} = \binom{m}{2} -
  \sum_{\kappa \in [k+1]} \binom{x_\kappa}{2}$,
  $$
  \alpha_{2m,\epsilon}^{\uparrow G}(q) =
  \alpha_{2m,\epsilon}^{\uparrow G}(q^{-1}) \; Y^{- \binom{m}{2} +
    \sum_{\kappa \in [k+1]} \binom{x_\kappa}{2}} .
  $$
  We shall show below that 
  \begin{equation} \label{eq_28}
    \alpha_{2m,\epsilon}^{\uparrow G}(q^{-1}) \; Y^{\sum_{\kappa \in [k+1]}
    \binom{x_\kappa}{2}} q^{-m} = \epsilon \sum_{H \subseteq [m]}
    (-1)^{m + \|H\|} c_{G,H} \; \alpha_{2m,\epsilon}^{\uparrow H}(q^{-1})
  \end{equation}
  From these equations \eqref{eq_27} follows
  readily.
  
  It remains to prove \eqref{eq_28}. An easy
  computation gives
  \begin{align*}
    Y^{\sum_{\kappa \in [k+1]} \binom{x_\kappa}{2}} q^{-m} & =
    \epsilon \frac{q^{-2m} + \epsilon q^{-m}}{1 + \epsilon q^{-m}}
    Y^{\sum_{\kappa \in [k+1]}
      \binom{x_\kappa}{2}} \\
    & = \epsilon \frac{1}{1 + \epsilon q^{-m}} \left( Y^{\sum_{\kappa
          \in [k+1]} \binom{x_\kappa +1}{2}} + \epsilon
      q^{-m} Y^{\sum_{\kappa \in [k+1]} \binom{x_\kappa}{2}} \right) \\
    & = \epsilon \left( \frac{Y^{\sum_{\kappa \in [k+1]}
          \binom{x_\kappa + 1}{2}} - Y^{\sum_{\kappa \in [k+1]}
          \binom{x_\kappa}{2}}}{1 + \epsilon q^{-m}} + Y^{\sum_{\kappa
          \in [k+1]} \binom{x_\kappa}{2}} \right).
  \end{align*}
  From Proposition~\ref{proposition_3} we see that for
  $H \subseteq [m]$,
  \begin{equation*}
    \alpha_{2m,\epsilon}^{\uparrow H}(q^{-1}) =
    \begin{cases}
      \alpha_{2m+1}^{\uparrow H}(q^{-1}) & \text{if $N(H) = m$,} \\
      \alpha_{2m+1}^{\uparrow H}(q^{-1}) / (1 + \epsilon q^{-m}) &
      \text{otherwise.}
    \end{cases}
  \end{equation*}
  In view of \eqref{eq_23} and \eqref{eq_26}, we
  thus obtain
  \begin{align*}
    \alpha_{2m,\epsilon}^{\uparrow G} & (q^{-1}) Y^{\sum_{\kappa \in
        [k+1]} \binom{x_\kappa}{2}} q^{-m} \\
    & = \epsilon \alpha_{2m + 1}^{\uparrow G}(q^{-1}) \left(
      \frac{Y^{\sum_{\kappa \in [k+1]} \binom{x_\kappa + 1}{2}} -
        Y^{\sum_{\kappa \in [k+1]} \binom{x_\kappa}{2}}}{1 + \epsilon
        q^{-m}} + Y^{\sum_{\kappa \in [k+1]} \binom{x_\kappa}{2}}
    \right) \\
    & = \epsilon \sum_{H \subseteq [m], N(H) \not = m} (-1)^{m +
      \|H\|} c_{G,H} \; \alpha_{2m+1}^{\uparrow
      H}(q^{-1}) / (1 + \epsilon q^{-m}) \\
    & \quad + \; \epsilon \sum_{H \subseteq [m], N(H) = m} (-1)^{m +
      \|H\|} c_{G,H} \; \alpha_{2m+1}^{\uparrow
      H}(q^{-1}) \\
    & = \epsilon \sum_{H \subseteq [m]} (-1)^{m + \|H\|} c_{G,H} \;
    \alpha_{2m,\epsilon}^{\uparrow H}(q^{-1}).
  \end{align*}
  This proves \eqref{eq_28}.
\end{proof}


\section{A conjecture for the orthogonal case}\label{section_conjecture}

In this section we discuss in more detail
Conjecture~\ref{conjecture_C}. As in Section~\ref{section_orthogonal},
let $\calV =(V,B,f)$ be an $n$-dimensional, non-degenerate quadratic
space over the finite field~$F = \F_q$. Our aim is to give, for $J
\subseteq [n-1]$, an expression for the
polynomial~$\alpha^J_{\calV}(q^{-1})$ in terms of parabolic length
functions on the Coxeter group~$W$ of type $A_{n-1}$.  If
Conjecture~\ref{conjecture_C} holds, the orthogonal case of
Theorem~\ref{theorem_A} follows directly from Theorem~\ref{theorem_1}.

Fix the Coxeter system $(W,S)$ where $W = \sym_n$ and $S = \{ s_1,
\dots, s_{n-1} \}$ denotes the standard set of Coxeter generators $s_i
= (i \;\; i+1 )$, $i \in [n-1]$. A crucial role is played by the
following statistic on $W$.

\begin{definition}[Length $L$]
  Recalling the notation from Section~\ref{section_coxeter}, for
  $w \in W$ set
  \begin{equation}\label{eq_29}
    L(w) := \bfb \cdot \bfl_\upR(w), \quad \text{where }
    \bfb = \left(b_I \right)_{I \subseteq S} = \left((-1)^{|I|}
    2^{|S|-|I|-1} \right)_{I \subseteq S}. 
  \end{equation}
\end{definition}

It is well-known that the ordinary Coxeter length of a permutation $w
\in W$ is equal to the number of inversion pairs associated to $w$,
i.e.\ $l(w) = |\mathcal{I}(w)|$ where
$$
\mathcal{I}(w) := \left\{ (i,j) \mid 1 \leq i < j \leq n, i^w > j^w \right\}.
$$
The parabolic length function $L$ also has a simple interpretation in
terms of inversion pairs.

\begin{lemma} 
  For each $w \in W$,
  $$
  L(w) = | \{ (i,j) \in \mathcal{I}(w) \mid i \not \equiv j \mod 2 \} |.
  $$
\end{lemma}

\begin{proof}
  Let $w \in W$ and note that for any $I \subseteq [n-1]$,
  $$
  l^I_\upR(w) = |\left\{(i,j)\in\mathcal{I}(w) \mid [i,j-1] \not
    \subseteq I \right\}|.
  $$
  From this we derive
  \begin{align*}
    L(w) & = \frac{1}{2} \sum_{I\subseteq [n-1]} (-1)^{|I|}
    2^{|S|-|I|} l^I_\upR(w) \\
    & = \frac{1}{2} \sum_{(i,j) \in \mathcal{I}(w)} \sum_{I \subseteq
      [n-1]} (-1)^{|I|} 2^{|S|-|I|} \; \delta(\text{`} [i,j-1]
    \not\subseteq
    I \text{'}) \\
    & = \frac{1}{2} \sum_{(i,j) \in \mathcal{I}(w)} \left( \sum_{I
        \subseteq [n-1]} (-1)^{|I|} 2^{|S|-|I|} -
      \sum_{[i,j-1] \subseteq I} (-1)^{|I|} 2^{|S|-|I|} \right) \\
    & = \frac{1}{2} \sum_{(i,j) \in \mathcal{I}(w)} \left( (2-1)^{|S|}
      -
      (-1)^{j-i} (2-1)^{|S|-(j-i)} \right) \\
    & = \frac{1}{2} \sum_{(i,j) \in \mathcal{I}(w)}
    \left(1-(-1)^{j-i}\right),
  \end{align*}
  where the Kronecker-delta $\delta(\text{`} [i,j-1] \not\subseteq I
  \text{'}) \in \{1,0\}$ reflects whether or not the inclusion
  $[i,j-1] \not\subseteq I$ holds.
\end{proof}

\begin{definition}[Chessboard elements]
  We say that $w \in W$ is a \emph{chessboard element} if $i + i^w
  \equiv j + j^w$ modulo $2$ for all $i,j \in [n]$. Clearly, the set
  $\mathcal{C}_n$ of chessboard elements forms a subgroup of $W$.
  Note that $\mathcal{C}_n$ contains a subgroup~$\mathcal{C}_{n,0}$
  consisting of elements~$w$ such that $i \equiv i^w$ modulo $2$ for
  all $i \in [n]$. If $n = 2m+1$ is odd, we have $\mathcal{C}_n =
  \mathcal{C}_{n,0} \cong \sym_{m + 1} \times \sym_m$. If $n = 2m$ is
  even, we have $\mathcal{C}_n = \langle w_0 \rangle \ltimes
  \mathcal{C}_{n,0}$, where $w_0$ denotes the longest element of $W$,
  and $\mathcal{C}_{n,0} \cong \sym_m \times \sym_m$.
  
  We write~$\sigma: W \rightarrow \{1,-1\}$, $w \mapsto (-1)^{l(w)}$
  for the \emph{sign character}, and $\tau : \mathcal{C}_n \rightarrow
  \{1,-1\}$ for the linear character with $\ker(\tau) =
  \mathcal{C}_{n,0}$.  Recall from the introduction that, in the
  even-dimensional case, we attach a sign~$\epsilon \in \{1,-1\}$
  to~$\calV$. Observing that $\tau$ is trivial for $n$ odd, we define
  $$
  \chi_\epsilon : \mathcal{C}_n \rightarrow \{1,-1\}, \chi_\epsilon (w) :=
  \begin{cases}
    \sigma(w) & \text{if $n$ is odd, or if $n$ is even and $\epsilon = 1$,} \\
    \sigma(w) \tau(w) & \text{if $n$ is odd, or if $n$ is even and
      $\epsilon = -1$.}
  \end{cases}
  $$
\end{definition}

\addtocounter{thmx}{-1} 

\begin{conjecture}
  For each $J \subseteq [n-1]$,
  \begin{equation} \label{eq_30}
    \alpha^J_{\calV}(q^{-1}) = \alpha^J_{n,\epsilon}(q^{-1}) =
        \sum_{\substack{w \in \mathcal{C}_n \\ 
        D_\upL(w) \subseteq J}} \chi_\epsilon (w) q^{-L(w)} .
  \end{equation}
\end{conjecture}

Note that, if Conjecture~\ref{conjecture_C} holds, the orthogonal case
of Theorem~\ref{theorem_A} follows from Theorem~\ref{theorem_1},
equation~\eqref{eq_5}, with~$W' = \mathcal{C}_n$, $\bfb$ as defined
in~\eqref{eq_29} and $\chi = \chi_\epsilon$.
Conjecture~\ref{conjecture_C} has been confirmed for $|J| \leq 1$ and
verified for~$n \leq 13$.

\begin{example} 
  For $n=3$, we have
  \begin{figure}[H]
  \begin{center}
  \begin{tabular}{|l!{\vrule width 1pt}c|c|c|c|}
    \hline
    $J \subseteq [2]$&$\varnothing$&$\{1\}$&$\{2\}$&$\{1,2\}$\\
    \hline\hline
    $a^J_3(q)$&$1$&$q^2$&$q^2$&$q^3-q$\\
    \hline
    $\alpha^J_3(q^{-1})$&$1$&$1$&$1$&$1-q^{-2}$\\
    \hline
  \end{tabular}
  \end{center}
  \end{figure}
  For $w \in W = \langle s_1, s_2 \rangle = \sym_3$, the statistic
  $L(w) = 2l(w) - l_{\upR}^{\{1\}}(w) - l_{\upR}^{\{2\}}(w)$, the
  character $\chi_\epsilon (w) = \sigma(w) = (-1)^{l(w)}$ and the left
  descent set $D_\upL(w)$ take the values
  \begin{figure}[H]
  \begin{center}
  \begin{tabular}{|l!{\vrule width 1pt}c|c|c|c|c|c|}
  \hline
  $w$&$\Id$&$s_1$&$s_2$&$s_1s_2$&$s_2s_1$&$s_1s_2s_1$\\
  \hline\hline
  $L(w)$&$0$&$1$&$1$&$1$&$1$&$2$\\
  \hline
  $\chi_\epsilon(w)$&$1$&$-1$&$-1$&$1$&$1$&$-1$\\
  \hline
  $D_\upL(w)$&$\varnothing$&$\{1\}$&$\{2\}$&$\{1\}$&$\{2\}$&$\{1,2\}$\\
  \hline
  \end{tabular}
  \end{center}
  \end{figure}
\end{example}

If $n$ is odd or if $n$ is even and $\epsilon = 1$, the character
$\chi_\epsilon$ naturally extends to the sign character on the whole
group $W$. Interestingly, in this case also a modified version of
equation~\eqref{eq_30} seems to hold, where $\chi_\epsilon$ is
replaced by $\sigma$ and one sums over \emph{all} elements $w \in W$.
In fact, we originally introduced chessboard elements in an attempt to
control cancellation in this larger sum.  Evidently, the contributions
of any two elements $w_1, w_2 \in W$ with $w_1^{-1} w_2 \in S$ and
$L(w_1) = L(w_2)$ cancel each other.  Therefore we were led to sum
over the set
$$
\mathcal{M} := \{w \in W \mid \forall s\in S: D_\upL(w) \not = D_\upL(ws)
\text{ or } L(w)\not = L(ws) \}.
$$
The set~$\mathcal{M}$ is easily seen to be closed under
right-multiplication by the longest element $w_0$ and might indeed
coincide with~$\mathcal{C}_n$. Aided by computer evidence, we
distilled Conjecture~\ref{conjecture_C} out of this circle of ideas.


\appendix

\section{Explicit examples}\label{appendix}

Here we collect a few examples of Igusa-type functions for the family
$$
\bfF = \left(F_J(\bfX)\right)_{J\subseteq[n-1]} = \left(\prod_{i\in
    J}\frac{X_i}{1-X_i}\right)_{J\subseteq[n-1]}.
$$


\subsection{A flag of alternating bilinear forms} 
Let $V$ be a vector space of dimension $n = 6$, equipped with a flag
of alternating bilinear forms $\bfB$ of type $I = \{4\}$. This kind of
example was considered immediately after
Definition~\ref{definition_4}. We have

\begin{figure}[H]
\begin{center}
\begin{tabular}{|l!{\vrule width 1pt}c|c|c|c|}
 \hline
 $J \subseteq [4]$&$\varnothing$&$\{2\}$&$\{4\}$&$\{2,4\}$\\
 \hline\hline
 $a^J_{6,\{4\}}(q)$&$1$&$q^8+q^4+q^2$&$q^8+q^6+1$&
 $q^2 (q^{12} - 1)/(q^2 -1)$\\
 \hline 
 $\alpha^J_{6,\{4\}}(q^{-1})$ & $1$ & $1 + q^{-4} + q^{-6}$ & $1 +
 q^{-2} + q^{-8}$ & $(1 - q^{-12})/(1 - q^{-2})$ \\
 \hline 
\end{tabular}
\end{center}
\end{figure}

Thus, in view of Definition~\ref{definition_5},
$$
\Ig_{(V,\bfB)}(q^{-1},\bfX)=\frac{1 +(q^{-4}+q^{-6}) X_1 +
  (q^{-2}+q^{-8})X_2 + q^{-10}X_1X_2}{(1-X_1)(1-X_2)}.
$$
Note that, in accordance with Theorem~\ref{theorem_B},
\begin{equation*}
\Ig_{(V,\bfB)}(q,X_1^{-1},X_2^{-1}) = q^{10} \,
\Ig_{(V,\bfB)}(q^{-1},X_2,X_1).  
\end{equation*}


\subsection{Quadratic spaces}

\subsubsection{Odd dimension}

Let $V$ be a vector space of odd dimension $n = 3$, equipped with a
non-degenerate quadratic form $f$. We have

\begin{figure}[H]
\begin{center}
\begin{tabular}{|l!{\vrule width 1pt}c|c|c|}
 \hline
 $J \subseteq[2]$&$\varnothing$&$\{1\}, \{2\}$&$\{1,2\}$\\
 \hline\hline
 $a^J_3(q)$&$1$&$q^2$&$q^3-q$\\
 \hline
 $\alpha^J_3(q^{-1})$&$1$&$1$&$1-q^{-2}$\\
 \hline
\end{tabular}
\end{center}
\end{figure}
In view of Definition~\ref{definition_6}, 
$$
\Ig_3(q^{-1},\bfX) = \frac{1-q^{-2}X_1X_2}{(1-X_1)(1-X_2)}
$$
and this satisfies the functional equation predicted by
Theorem~\ref{theorem_A}.


\subsubsection{Even dimension}

Let $V$ be a vector space of dimension $n = 4$, equipped with a
non-degenerate quadratic form $f$. We have

\begin{figure}[H]
\begin{center}
\begin{tabular}{|l!{\vrule width 1pt}c|c|c|c|c|}
 \hline
 $J \subseteq [3]$ & $\varnothing$ & $\{1\},\{3\}$ & $\{2\}$ &
 $\{1,2\},\{2,3\}, \{1,3\}$ & $\{1,2,3\}$ \\
 \hline\hline 
 $a^J_{4,1}(q)$ & $1$ & $q^3-q$ & $q^4+q^2$ & $q^5-q^3$ &
 $q^2(q^2-1)^2$ \\
 \hline
 $\alpha^J_{4,1}(q^{-1})$ & $1$ & $1-q^{-2}$ & $1+q^{-2}$ & $1-q^{-2}$
 & $(1-q^{-2})^2$ \\
 \hline
 $a^J_{4,-1}(q)$ & $1$ & $q^3+q$ & $q^4+q^2$ & $q^5+q^3$ &
 $q^2(q^4-1)$ \\
 \hline
 $\alpha^J_{4,-1}(q^{-1})$ & $1$ & $1+q^{-2}$ & $1+q^{-2}$ &
 $1+q^{-2}$ & $1-q^{-4}$ \\
 \hline
\end{tabular}
\end{center}
\end{figure}

Note that, for both values of~$\epsilon$, the number of
\emph{distinct} polynomials among the $\alpha^J_{4,\epsilon}(q^{-1})$
is comparatively small.  This illustrates that the map $J \mapsto
\alpha^J_{n,\epsilon}(q^{-1})$ factors over the bisecting map $\phi$;
cf.~Proposition~\ref{proposition_3}. One readily
computes
{\small
\begin{align*}
  \Ig_{4,1}(q^{-1},\bfX) & = \frac{1 + q^{-2} \left(-X_1X_2 + X_1X_3 -
      X_2X_3 - X_1 + X_2 - X_3 \right) + q^{-4} X_1X_2X_3}{(1 - X_1)(1
    - X_2)(1 - X_3)}, \\
  \Ig_{4,-1}(q^{-1},\bfX) & = \frac{1 +
    q^{-2}\left(-X_1X_2-X_1X_3-X_2X_3 +X_1+X_2+X_3\right) -
    q^{-4}X_1X_2X_3}{(1-X_1)(1-X_2)(1-X_3)}.
\end{align*}
}
These Igusa-type functions satisfy the functional equations predicted
by Theorem~\ref{theorem_A}.


\subsection{A symmetric bilinear space in characteristic $2$}

Let $V = \langle e_1, e_2,e_3,e_4 \rangle$ be a $4$-dimensional vector
space over a finite field $F = \F_q$ with $\cha F = 2$. Let $B$ denote
the non-degenerate symmetric bilinear form on $V$ such that
$B(e_i,e_j) = \delta_{ij}$, and define $f : V \rightarrow F$ by $f(x)
:= B(x,x)$. We consider flags of non-degenerate symmetric bilinear
subspaces in the non-degenerate symmetric bilinear space $\calV :=
(V,B,f)$. For $J \subseteq \{1,2,3\}$ let $a_\calV^J(q)$ denote the
number of non-degenerate flags of type $J$ in $\calV$, and define the
normalised polynomials $\alpha_\calV^J(q^{-1})$ and the `Igusa-type
function' $\Ig_\calV(q^{-1},\bfX)$ accordingly. Then we have

\begin{figure}[H]
\begin{center}
\begin{tabular}{|l!{\vrule width 1pt}c|c|c|c|c|}
 \hline
 $J \subseteq [3]$ & $\varnothing$ & $\{1\},\{3\}$ & $\{2\}$ &
 $\{1,2\},\{2,3\}, \{1,3\}$ & $\{1,2,3\}$ \\
 \hline\hline 
 $a^J_\calV(q)$ & $1$ & $q^3$ & $q^4+q^2$ & $q^5$ & $q^4(q^2-1)$ \\
 \hline
 $\alpha^J_\calV(q^{-1})$ & $1$ & $1$ & $1+q^{-2}$ & $1$ & $1-q^{-2}$ \\
 \hline
\end{tabular}
\end{center}
\end{figure}
The associated `Igusa-type function' 
$$
\Ig_\calV(q^{-1},\bfX) = \frac{1 + q^{-2} X_2 (1 - X_1 -
  X_3)}{(1-X_1)(1-X_2)(1-X_3)}
$$
\emph{does not} satisfy a functional equation; this illustrates our
remarks at the end of the Introduction.


\bibliographystyle{amsplain}

\providecommand{\bysame}{\leavevmode\hbox to3em{\hrulefill}\thinspace}
\providecommand{\MR}{\relax\ifhmode\unskip\space\fi MR }
\providecommand{\MRhref}[2]{%
  \href{http://www.ams.org/mathscinet-getitem?mr=#1}{#2}
}
\providecommand{\href}[2]{#2}

\end{document}